\documentclass{amsart}
\usepackage{graphicx}
\usepackage{amssymb}
\usepackage{epstopdf}
\newtheorem{thm}{THEOREM}[section]

\newtheorem{cor}[thm]{COROLLARY}
\newtheorem{defn}[thm]{DEFINITION}
\newtheorem{ex}[thm]{EXAMPLE}

\newtheorem{lemma}[thm]{LEMMA}

\newtheorem{prop}[thm]{PROPOSITION}

\newcommand{\bA}{{\mathbf  A}}
\newcommand{\bH}{{\mathbf  H}}
\newcommand{\bO}{{\mathbf  O}}
\newcommand{\bS}{{\mathbf  S}}
\newcommand{\bU}{{\mathbf  U}}
\newcommand{\cA}{{\mathcal A}}
\newcommand{\cat}{\mbox{\rm cat}}

\newcommand{\catG}{\mbox{\rm cat}_{G}}
\newcommand{\catO}{\mbox{\rm cat}_{\Oq}}
\newcommand{\catt}{\cat_{\subtrans}}
\newcommand{\cH}{{\mathcal H}}
\newcommand{\cI}{{\mathcal I}}
\newcommand{\cK}{{\mathcal K}}
\newcommand{\cN}{{\mathcal N}}
\newcommand{\cP}{{\mathcal P}}

\newcommand{\cT}{{\mathcal T}}
\newcommand{\cU}{{\mathcal U}}
\newcommand{\cV}{{\mathcal V}}
\newcommand{\cW}{{\mathcal W}}
\newcommand{\cX}{{\mathcal X}}
\newcommand{\cZ}{{\mathcal Z}}
\newcommand{\ds}{\displaystyle}
\newcommand{\E}{{\mathcal E}}
\newcommand{\ecatt}{\cat^e_{\subtrans}}
\newcommand{\eop}{{\;\;\;  \Box}}
\newcommand{\F}{{\mathcal F}}
\newcommand{\Fr}{\mathbf{Fr}}
\newcommand{\G}{\Gamma}
\newcommand{\Isom}{\mathbf{Isom}}
\newcommand{\mD}{{\mathbb D}}
\newcommand{\mF}{{\mathbb F}}
\newcommand{\mQ}{{\mathbb Q}}
\newcommand{\mR}{{\mathbb R}}
\newcommand{\mS}{{\mathbb S}}
\newcommand{\mT}{{\mathbb T}}
\newcommand{\mZ}{{\mathbb Z}}
\newcommand{\oF}{\overline{{\mathcal F}}}
\newcommand{\oL}{\overline{L}}
\newcommand{\Oq}{{\mathbf O(q)}}
\newcommand{\owhL}{\overline{\widehat L}}
\newcommand{\sop}{{\bf Proof:} }

\newcommand{\Up}{{\Upsilon}}

\newcommand{\what}{\widehat}

\newcommand{\whE}{{\widehat{\E}}}
\newcommand{\whF}{{\widehat{\F}}}

\newcommand{\whg}{{\widehat{g}}}
\newcommand{\whH}{{\widehat H}}
\newcommand{\whL}{{\widehat L}}
\newcommand{\whM}{{\widehat M}}
\newcommand{\whpi}{{\widehat{\pi}}}
\newcommand{\whQ}{{\widehat{Q}}}
\newcommand{\whU}{{\widehat U}}
\newcommand{\whUp}{{\widehat{\Upsilon}}}
\newcommand{\whV}{{\widehat V}}
\newcommand{\whW}{{\widehat W}}
\newcommand{\whw}{{\widehat{w}}}
\newcommand{\whx}{{\widehat{x}}}
\newcommand{\why}{{\widehat{y}}}
\newcommand{\whz}{{\widehat{z}}}
\newcommand{\whZ}{{\widehat{Z}}}
\newcommand{\wtB}{{\widetilde B}}
\newcommand{\wtM}{{\widetilde M}}
\newcommand{\wtpi}{{\widetilde{\pi}}}
\newcommand{\wtW}{{\widetilde W}}
\def\subtrans{\mathbin{\cap{\mkern-9mu}\mid}\,\,}

\begin{document}

\title{Transverse LS-category for Riemannian foliations}

\author{Steven Hurder}
\thanks{SH supported in part by   NSF grant  DMS-0406254}
\address{Steven Hurder, Department of Mathematics, University of Illinois at Chicago, 322 SEO (m/c 249), 851 S. Morgan Street, Chicago, IL 60607-7045, USA}
\email{hurder@uic.edu}
\author{Dirk T\"oben}
 \thanks{DT supported  by the Schwerpunktprogramm SPP 1154  of the DFG}
\address{Dirk T\"oben, Mathematisches Institut, Universit\"at zu K\"oln, Weyerthal 86-90, 50931 K\"oln, Germany}
\email{dtoeben@math.uni-koeln.de}
\thanks{Preprint date: October 1, 2006}

\date{}



\begin{abstract}
We study the transverse  Lusternik-Schnirelmann category of a Riemannian foliation $\F$ on a compact manifold $M$. We obtain necessary and sufficient conditions for when the transverse  category $\catt(M,\F)$ is finite. We also introduce a variation on the concept of transverse LS category, the essential transverse  category $\ecatt(M,\F)$, and show that this is finite for every Riemannian foliation. Also, $\ecatt(M,\F) = \catt(M,\F)$   if  $\catt(M,\F)$    is finite. A generalization of the Lusternik-Schnirelmann theorem is also given: the essential transverse  category $\ecatt(M,\F)$  is  a lower bound for the number of critical leaf closures of a basic $C^1$-function on $M$.   
\end{abstract}
 
 \maketitle

\tableofcontents

 \vfill
\eject

\section{Introduction} \label{sec-intro}

Let $f \, \colon M \to \mR$ be a $C^1$-function on a closed Riemannian manifold $M$. We say that $x \in M$ is a critical point if the gradient $\nabla f$ vanishes at $x$. A well-known formula of  Lusternik-Schnirelmann 
\cite{LS1934,James1978,James1995,CLOT2003} gives a lower-bound estimate for the number of critical points, 
\begin{equation}\label{eq-LS}
\# \{ x  \mid x \in M  ~{\rm is ~ critical ~ for } ~ f\} \geq \cat(M)
\end{equation}
where $k = \cat(M)$ is the Lusternik-Schnirelmann category of $M$, which is defined    as   the least number of     open sets $\{U_1, \ldots , U_k\}$ required to cover $M$ such that each $U_{\ell}$ is contractible in $M$ to a point. The LS category  is    a measure of the topological complexity of $M$; in the case where $f$ is $C^2$ with only Morse-type singularities,  the value  $\cat(M)$  gives a lower bound for the number of cells in a cellular decomposition.     

Let $G$ be a compact Lie group and suppose there is a smooth action $M \times G \to M$ which we can assume preserves a Riemannian metric on $M$. A $C^1$-function $f \, \colon M \to \mR$ is $G$-invariant if $f(x \, A) = f(x)$ for all $A \in G$, hence the set of critical points $\nabla f = 0$ is $G$-invariant. 
Each $G$-orbit $x\, G$  is a closed submanifold of $M$, and the number of critical $G$-orbits is estimated by the $G$-category, 
\begin{equation}\label{eq-GLS}
\# \{ x \, G \mid x \in M ~ {\rm is ~ critical ~ for } ~ f\} \geq \cat_G(M)
\end{equation}
where now $\cat_G(M)$ is the   least number of    $G$-invariant  open sets $\{U_1, \ldots , U_k\}$ required to cover $M$ such that each $U_{\ell}$ is 
$G$-contractible in $M$ to a single orbit  (see 
Marzantowicz \cite{Marzantowicz1989} for actions of a compact Lie group,  and Ayala, Lasheras and Quintero  \cite{ALQ2001} for   proper actions of   Lie groups.)

 A smooth action of a non-compact Lie group $G$ on a compact manifold is never proper, and the study of its LS-category theory in this case is much more difficult,   as the analysis must take into account the dynamical properties of the action.   
In order that  the condition that the critical sets  $\nabla f = 0$ for the gradient  be $G$-invariant, it is useful to assume the action of the group $G$ preserves a Riemannian metric on $M$, and thus the orbits of the action define a 
\emph{singular Riemannian foliation} (SRF) of $M$ 
\cite{Haefliger1989,Molino1988,Molino1994}. If all orbits of the action have the same dimension, then the orbits define a Riemannian foliation of $M$.

 In this paper, we study the LS-category theory for Riemannian foliations, and  prove  a Lusternik-Schnirelmann type  estimate   for the case of Riemannian foliations. 
 \begin{thm} \label{thm-LSR} 
 Let $\F$ be a Riemannian foliation for a compact manifold $M$, and 
  let  $f \, \colon M \to \mR$ be a $C^1$-map which is constant along the leaves of $\F$.  A leaf $L_x$ of $\F$ through a point $x \in M$  is critical if $\nabla f | L_x$ = 0, and hence $\nabla f$ also  vanishes on the leaf closure $\oL_x$. Then the number of critical leaf closures has a lower bound estimate    by the 
  \emph{essential transverse LS category} $ \ecatt(M,\F)$ of $\F$, 
\begin{equation}\label{eq-FLS}
\# \{ \oL_x  \mid x \in M  ~{\rm is ~ critical ~ for } ~ f\} \geq \ecatt(M,\F)
\end{equation}
 \end{thm}

 In the case where all leaves of $\F$ are compact,   Colman      \cite{Colman1998,CM2001}  proved  a   lower bound estimate for the number of critical leaves in terms of the transverse LS category $\catt(M,\F)$, which the estimate (\ref{eq-FLS}) generalizes to the general case when $\F$ has non-compact leaves. 
 The transverse LS category $\catt(M,\F)$ is   infinite when $\F$  has no compact leaves, while the essential transverse LS category  $\ecatt(M,\F)$ introduced in this paper is always   finite.   
 
 \eject

The  basic concept is that of a {\it foliated homotopy}\,: 
given foliated manifolds  $(M,{\F})$ and $(M',{\F}')$, a  map $f \, \colon M \to M'$ is 
  said to be {\it  foliated} 
if for each leaf $L \subset M$ of $\F$, there exists a leaf $L' \subset M'$ of $\F'$ such that $f(L) \subset L'$. A
 $C^r$-map $H \, \colon M' \times [0,1] \to M$, for $r \geq 0$, is said to be a \emph{foliated $C^r$-homotopy}
if  $H_t$  is foliated for    all  $0 \leq t \leq 1$,  and $H_0(x) = x$ for all $x \in U$. As usual,  $H_t(x) = H(x,t)$.

Unless otherwise specified, we assume that all maps and homotopies are smooth.

Let  $U \subset M$ be an open saturated subset. We say that $U$ is {\it transversely categorical} if  there 
is a foliated  homotopy $H \, \colon U \times [0,1] \to M$ such that 
  $H_1\, \colon U \to M$ has image in a 
single leaf of $\F$.  

\begin{defn}\label{def-cat}
The  transverse  LS category $\catt  (M,\F)$ of a foliated  manifold 
$(M,\F)$ is the least number of transversely categorical open saturated sets 
required to cover $M$. If no such covering exists, then  set
$\catt  (M,\F)=\infty$.
\end{defn}

The basic properties of transverse LS category are given in  \cite{Colman1998,CM2001}. If a foliation $\F$ is defined by a fibration $M \to B$ over a compact manifold $B$, then $\catt  (M,\F)= \textrm{cat}(B) < \infty$, so the LS category of $\F$ agrees with the LS category of the leaf space $M/\F$ in this case. Also, the transverse LS category is an invariant of foliated homotopy. 
The transverse saturated category $\catt(M,\F)$   has been  further   studied by various authors 
 \cite{Colman2002b,Colman2004,CH2004,HT2006b,HW2006,LW2002,Wolak2002}.

The assumption that $\catt  (M,\F)$ is finite is a strong hypothesis on $\F$, and has  consequences for the dynamical properties of $\F$. 
Let  $\{U_1, \ldots , U_k\}$ be  a minimal  cardinality  covering of $M$ by categorical open $\F$-saturated sets, so that  $k = \catt(M,\F)$. Then each $U_i$ contains a compact minimal set $K_i \subset U_i$ for $\F$. The first author  showed in \cite{Hurder2006a}  that if $H_{\ell} \, \colon U_{\ell} \times [0,1] \to M$ is a foliated homotopy, 
the image $H_{\ell, t}(K_{\ell})$ is a compact minimal set for all $t$, and in particular $H_{\ell , 1}(K_{\ell})$ must be a compact minimal set. Thus, if $H_{\ell}$ is a categorical homotopy, then the image $H_{\ell , 1}$ must be contained in a compact leaf of $L$. If $\F$ has no compact leaves, or simply not enough compact leaves, then   a categorical covering of $M$ cannot be found. 
 
Let  $U \subset M$ be an open saturated subset. We say that $U$  is {\it essentially transversely categorical} if  there 
is a foliated  homotopy $H \, \colon U \times [0,1] \to M$ such that 
  $H_1\, \colon U \to M$ has image in a 
minimal set  of $\F$.  
\begin{defn}\label{def-ecat}
The  \emph{essential transverse  LS category} $\ecatt  (M,\F)$ of a foliated  manifold 
$(M,\F)$ is the least number of essentially transversely categorical open saturated sets 
required to cover $M$. If no such covering exists, then   set 
$\ecatt  (M,\F)=\infty$.
\end{defn}

 With this definition, we obtain the following fundamental result:
 
  \begin{thm}\label{thm-main1} Let $\F$ be a Riemannian foliation of a compact smooth manifold $M$. 
  Then the essential transverse       category $\ecatt (M,\F)$ is finite. 
  If the transverse category  $\catt  (M,\F)$ is   finite, then   $\catt  (M,\F) =   \ecatt  (M,\F) $.
   \end{thm} 
   
   We obtain an exact characterization of which Riemannian foliations have  $\catt  (M,\F)$. 
   Let  $L$  be  a leaf of $\F$. A  \emph{foliated isotopy} of $L$ is a smooth map  $I \, \colon L \times [0,1] \to M$ such that $I_0 \colon L \to M$ is the inclusion of $L$, and for each $0 \leq t \leq 1$, $I_t \colon L \to M$ is a diffeomorphism onto its image $L_t$, which is a leaf of $\F$. We   say that the image leaf $L_1$ is \emph{foliated isotopic} to $L$.     
   
    Let $\cI_L$ denote the set of leaves of $\F$ which are foliated isotopic to $L$. For $x \in M$, we set $\cI_x = \cI_{L_x}$. 
    
    For an arbitrary foliation, one cannot expect the isotopy classes $\cI_x$ to have any nice properties at all, and   typically one expects that $\cI_x = L_x$.  However, for a Riemannian foliation, each isotopy class $\cI_x$ is a smooth submanifold of $M$, and 
     the set of isotopy classes of the leaves of $\F$ defines a  Whitney stratification of $M$. (This is proven in  section~\ref{sec-strata}.)   A leaf $L_x$ (as well as its corresponding stratum $\cI_x$) is said to be \emph{locally minimal} if $\cI_x$ is a \emph{closed} submanifold of $M$.    
 
  \begin{thm}\label{thm-main2} Let $\F$ be a Riemannian foliation of a compact smooth manifold $M$. 
  Then $\catt (M,\F)$ is finite if and only if each locally minimal leaf $L_x$ is compact, and hence  the locally minimal set  $\cI_x$  is a union of compact leaves.
     \end{thm} 
   
 The proofs of Theorems~\ref{thm-LSR}, \ref{thm-main1} and \ref{thm-main2} are based on the very special geometric properties of Riemannian foliations, which were developed by Molino in a series of papers  \cite{Molino1977,Molino1982,Molino1988,Molino1994}, and see also Haefliger \cite{Haefliger1985,Haefliger1988, Haefliger1989}. We first recall several key facts  in order to state our next result; details are given  in   section~\ref{sec-molino}.

One of the remarkable corollaries  of the Molino structure theory  is  that for a leaf $L$ of a Riemannian foliation, its closure $\oL$ is a minimal set for $\F$. Thus, the condition in Theorem~\ref{thm-LSR} that the critical points of a leafwise constant $C^1$-function consists of unions of leaf closures, means that in fact    they are unions of minimal sets.

 Molino's analysis of the geometry and structure of a Riemannian foliation is based on  the desingularization of $\F$ using the geometry of a  foliated   $\Oq$-bundle:     
Let $\whM \to M$ denote the principle $\Oq$-bundle of orthonormal frames for the normal bundle to $\F$.
  There exists an $\Oq$-invariant  foliation $\whF$ on $\whM$ whose leaves are the holonomy coverings of the leaves of $\F$, and in particular have the same dimension as the leaves of $\F$.  The closures of the leaves of $\whF$ are the fibers of an $\Oq$-equivariant fibration $\whUp \, \colon \whM \to \whW$, where $\whW$ is the quotient space by the leaf closures of $\whF$, and the $\Oq$-action on $\whM$ induces the smooth action of $\whW$. The quotient  space 
$$\whW/\Oq ~ \cong ~ W \equiv  M/\oF$$
is naturally identified with the singular quotient space $W$ of $M$ by the leaf closures of $\F$. 
Given   $w \in \whW$, the inverse image $\widehat{\pi}^{-1}(w) = \owhL$   is the   closure of each leaf 
$\whL_{\whx} \subset \owhL$, and the projection of such $\owhL$ to $M$  is the closure $\oL$ of a leaf $L_x$ of $\F$.

The smooth  action of $\Oq$ on $\whW$ defines the orbit-type stratification of $\whW$ in terms of the stabilizer groups of the action, which is one of the key concepts for the study  of smooth actions of $\Oq$ 
(see \cite{Bredon1972, Davis1978, DK2000, HsHs1967, Janich1968, tomDieck1987}.)    A key idea for this work is the use of   the associated   Whitney stratification, 
 $\whW = \cZ_1 \cup \cdots \cup \cZ_K$, where each set $\cZ_{\ell}$ is a closed submanifold of $\whW$ which is $\Oq$-invariant and such that $\cZ_{\ell}/\Oq$ is connected (see 
sections~ \ref{sec-orbit}, \ref{sec-strata} and \ref{sec-finite}.) A stratum $\cZ_{\ell}$ is said to be 
\emph{locally minimal} if it is a closed submanifold.

Let $\catO(\whW)$ denote the  $\Oq$-equivariant category of the space $\whW$, which  is finite as $\whW$ is compact. We then have the following interpretation of the essential transverse category:
\begin{thm}\label{thm-main3}  Let $\F$ be   a Riemannian foliation of a compact manifold $M$. Then
\begin{equation}\label{eq-equicat}
\ecatt(M,\F) = \catO(\whW)
\end{equation}
\end{thm}
 
Note that for a smooth action of a compact group $G$ on a compact manifold $N$, the   $G$-equivariant category of $N$ defined using smooth $G$-homotopies is equal to the $G$-category defined using continuous homotopies,   as a continuous homotopy can be approximated by a smooth homotopy (Theorem~4.2, Chapter~VI of \cite{Bredon1972}.) Thus, 
 the  calculation of the    equivariant category $\catO(\whW)$ is a purely topological problem.

 The proof of Theorem~\ref{thm-main3} introduces one of the main new technical ideas of this paper, which can be called the ``synchronous lifting property''.
 The choice of a  projectable Riemannian metric on $M$, which is $\F$-projectable when restricted to the normal bundle to $\F$, defines a Riemannian metric on the principle $\Oq$-frame bundle $\pi \, \colon \whM \to M$ which is projectable with respect to the lifted foliation $\whF$.
 This projectable  Riemannian  metric   in turn defines a horizontal distribution in $T\whM$ which is transverse to the fibers of $\pi$. The first key idea is that a homotopy $H \colon U \times [0,1] \to M$ on $M$ can be lifted to an $\Oq$-equivariant $\whF$-foliated homotopy of $\whU = \pi^{-1}(U)$, $\whH  \, \colon \whU \times [0,1] \to \whM$. This shows that $\F$-categorical open sets on $M$ are equivalent to $\whF$-categorical, $\Oq$-equivariant categorical sets on $\whM$. 
 
 Similarly, the projection 
 $\whUp \, \colon \whM \to \whW$ from the frame bundle to the space of leaf closures is a Riemannian submersion, so has a natural horizontal distribution which is transverse to the projection $\whUp$.  This is used to show that  $\whF$-categorical, $\Oq$-equivariant categorical sets on $\whM$ are equivalent to $\Oq$-equivariant categorical sets on $\whW$.  
 
  It is interesting to note that our technique for   lifting homotopies to equivariant homotopies, used in both sections~\ref{sec-equivariant1} and \ref{sec-equivariant2},  is similar  to the method of ``averaging isotopies'' used in the proof of (Theorem~3.1, Chapter~VI of \cite{Bredon1972}).

The remarkable    aspect of these arguments is that the connection data is used to define equivariant lifts of the given homotopy; the     horizontal distributions are   used to ``synchronize'' the orthonormal frames along the traces of the homotopies.  This is possible, even  though the homotopy itself need not transform  normal frames to normal frames.  This technique  also     has applications to the theory of   secondary characteristic classes of $\F$, especially residue theory \cite{HT2006c,LP1976a,LP1976b}.

  This paper can be viewed as a sequel to the work \cite{CH2004} by the first author with Hellen Colman. The results of this paper also extends the results of \cite{Colman2004}.

 The research  collaboration of the authors resulting in this work was made possible by the generous support of the second author by the 
 Schwerpunktprogramm SPP 1154 ``Globale Differentialgeometrie'' of the Deutschen Forschungsgemeinschaft. Both authors are very grateful for this support.
 
 \vfill
\eject

\section{Transverse category}\label{sec-cat}

We assume that  $M$ is a smooth, compact Riemannian manifold without boundary of dimension $m = p +q$, and   $\F$ is a smooth Riemannian foliation of dimension $p$ and codimension $q$.  Given $x \in M$ we will denote by $L_x$ the leaf of $\F$ containing $x$.   

Let $\E$ denote the singular Riemannian foliation (SRF) of $M$ defined by the closures of the leaves of $\F$. (See Molino \cite{Molino1988,Molino1994} for properties of $\E$.)
The tangential distribution $E = T\E$ is   integrable and  satisfies the regularity  conditions formulated by Stefan \cite{Stefan1974,Stefan1980}.  Note that all leaves of $\E$ are compact. 

The notion of foliated homotopy extends naturally to the case of singular foliations, so that one can define the transverse category $\catt(M,\E)$ of $\E$. 
 We   recall two   topological lemmas due to Colman \cite{Colman2004} which are used to relate the transverse categories of $\F$ and $\E$.

  \begin{lemma}\label{lem-foliated}
     Let $(M,\F)$ and $(M',\F')$ be two foliated manifolds and 
     $f\, \colon M\to M'$ be a foliated continuous map. Let $\E$ denote the partition of $M$ by the closures of the leaves of $\F$, and $\E'$ the corresponding partition of $M'$.  Then $f$ is also    
     $\E$-foliated.
 \end{lemma}
\sop
Let $L \subset M$ be a leaf of $\F$, and $L' \subset M'$ the leaf of $\F'$ such that $f(L) \subset L'$. Then 
$\ds   f(\oL) \subset \overline{f(L)} \subset \overline{L'}$. $\hfill \eop$
 
 \medskip
   
The second lemma is based on the special property of Riemannian foliations    that the closure of every leaf in a Riemannian foliation is a minimal set. 
 \begin{lemma}\label{lem-closures} 
 Let $\F$ be a Riemannian foliation of a compact manifold $M$. 
  Let $U\subset M$ be an $\F$-saturated open set and $L$ be a leaf of the 
  Riemannian foliation $\F$ 
  such that $L\subset U$. Then $\oL\subset U$.   
\end{lemma}
\sop
Let $L\subset U$, and suppose that the closure  $\oL$ is not a subset of $U$.   Then there exists a leaf 
$L' \subset \oL$ such that $L' \not\subset U$, and as $U$ is saturated, it follows that $L' \subset M - U$. The complement $M - U$ is a closed saturated set, so $\overline{L'} \subset M - U$. But $L' \subset \oL$ and $\oL$ is a minimal set implies  that $\overline{L'} = \oL$,  so $L \subset \overline{L'} \subset M - U$, which is a contradiction.     
$\hfill \eop$

\begin{prop}\label{prop-induced}
    Let $U\subset M$ be a saturated open set. If $H\, \colon U\times 
    [0,1] \to M$ is an $\F$-homotopy, then   $H$ is also an $\E$-homotopy.   If $H_1$ has 
    image in a single leaf $L\in\F$, or more generally in a minimal set $K$ of $\F$,  then $H_1$ has image in the leaf $K = \overline{L}$ of $\E$. 
\end{prop}
\sop The open set  $U$ is a $\E$-saturated set    by Lemma~\ref{lem-closures}. 
The map $H_t$ is $\E$-foliated by Lemma~\ref{lem-foliated}. Then  the map $H_1$ has image in the
  closure $K = \oL$ by Lemma~\ref{lem-foliated}.  $\hfill \eop$

 \begin{cor}\label{cor-induced}
Let $\F$ be a Riemannian foliation of a compact manifold $M$, then 
\begin{equation}\label{eq-compare1}
\catt(\E) \leq \ecatt(\F) \leq \catt(\F)
\end{equation}
\end{cor}

\section{Geometry of Riemannian foliations}\label{sec-molino}

The Molino structure theory  is a remarkable collection of results about the geometry and topology of a Riemannian foliation on a compact manifold. We recall some of the main results, and in doing so establish the notation which will be used in later sections.  The  reader should consult Molino \cite{Molino1982,Molino1988,Molino1994}, Haefliger \cite{Haefliger1985,Haefliger1988, Haefliger1989}, or  Moerdijk and Mr{\v{c}}un \cite{MM2003} for further details, noting  that our notation is an amalgam of those used by these authors.

Let $M$ denote a compact, connected smooth manifold without boundary, and $\F$ a smooth Riemannian foliation of codimension $q$, with  tangential distribution $T\F$. 

Let $g$ denote a   Riemannian metric on $TM$ which is \emph{projectable} with respect to $\F$. Identify the normal bundle  $Q$   with the orthogonal space $T\F^{\perp}$, and let $Q$ have the restricted Riemannian metric $g_Q =  g | Q$.  For a vector $X \in T_xM$ let $X^{\perp} \in Q_x$ denote its orthogonal  projection.

Given a leafwise path $\gamma$ between points $x, y$ on a leaf $L$, the transverse holonomy $h_{\gamma}$ along $\gamma$ induces a   linear transformation $dh_x[\gamma] \, \colon Q_x \rightarrow Q_y$. The   fact that the Riemannian metric $g$ on $TM$ is projectable  is equivalent to the fact  that    the linear holonomy transformation $dh_x[\gamma]$ is an isometry for all such paths.

Let  $\{E_1, \ldots , E_q\}$ be an orthogonal basis for $\mR^q$. Fix $x \in M$. 
An orthonormal  frame for $Q_x$ is an  isometric  isomorphism $e \, \colon \mR^q \to Q_x$.
Let $\Fr(Q_x)$ denote the space of orthogonal frames of $Q_x$. 
Given $e \in \Fr(Q_x)$  and  $A \in \Oq$ we obtain a new frame $R_A(e) = e   A = e \circ A$  where $A \, \colon \mR^q \to \mR^q$ is the map induced by matrix multiplication. 

 
 The group $\Isom(Q_x)$ of  isometries of $Q_x$ also acts naturally on $\Fr(Q_x)$: for $h_x \in \Isom(Q_x)$ and $e \in \Fr(Q_x)$, we define $h_x e = h_x \circ e \, \colon \mR^q \to Q_x$.

 Let $\pi \, \colon \whM \to M$ denote the bundle of orthonormal frames for $Q$, where the fiber over $x \in M$ is  $\pi^{-1}(x) = \Fr(Q_x)$. By the above remarks, $\whM$ is a principle $\Oq$-bundle. We use the notation $\whx = (x,e) \in \whM$ where $e \in \Fr(Q_x)$.
 
 There is a canonical $\mR^q$-valued 1-form, the {\it Solder 1-form} $\theta \, \colon T\whM \to \mR^q$, defined as follows: for $X \in T_{\whx}\whM$, then $d_{\whx}\pi(X) \in T_xM$ and set 
 $$\theta(X) = e^{-1}\left(d_{\whx}\pi(X)^{\perp}\right) \in \mR^q$$
 Note that for $A \in \Oq$, $R_A^*\theta = A^{-1} \circ \theta$.
 
Let  $\nabla$ denote the Levi-Civita connection  on   $Q \to M$ defined using the   Riemannian metric. Let $\omega \, \colon T\whM \to \mathfrak{o}(q)$ denoted the associated Maurer-Cartan 1-form, with values in the Lie algebra $\mathfrak{o}(q)$ of $\Oq$. Recall that  $\omega$ is $Ad$-related: for $A \in \Oq$  and $X \in T\whM$, we have
\begin{equation}
R_A^*(\omega) = Ad(A^{-1}) \circ \omega ~ ; ~ R_A^*(\omega)(X) = \omega(dR_A (X)) = Ad(A^{-1})(\omega(X))
\end{equation}
Let $\cH = \ker(\omega)  \subset T\whM$ denote the horizontal distribution for $\omega$.  Then   $\cH$ is invariant under the right action of $\Oq$, and for all $\whx \in \whM$, the differential $d_{\whx}\pi \, \colon \cH_{\whx} \to T_{x}M$ is an isomorphism.

Define an $\Oq$-invariant  Riemannian metric $\whg$ on $T\whM$, by requiring that the restriction $d\pi \, \colon \cH \to TM$ be an isometry,  and the fibers of $\pi$ are orthogonal to $\cH$. The metric restricted to the fibers is induced  by  the bi-invariant  metric on $\Oq$ which is  defined by the inner product    on $\mathfrak{o}(q)$, where $\langle A, B \rangle  = \frac{1}{2} \, Tr (A\, B)$ for matrices $A, B \in \mathfrak{o}(q)$.

The  metric on $TM$  is projectable implies the restriction  $\nabla^L$ to $Q | L \to L$ is a flat connection for each leaf $L \subset M$, so the horizontal distribution of $\nabla^L$ is integrable.
The inverse image $\pi^{-1}(L) \subset \whM$ is an $\Oq$-principle bundle over $L$, and the restricted  flat connection $\nabla^L$ defines an $\Oq$-invariant  foliation of  $\pi^{-1}(L)$,   whose leaves   cover $L$. The union of all such flat subbundles in $\whM$ defines an $\Oq$-invariant  foliation  $\whF$ of $\whM$ whose tangent distribution $T\whF$ is an integrable subbundle of $\cH$. The metric  on $T\whM$ is projectable for $\whF$, hence   $\whF$ is also a Riemannian foliation.

The direct sum of the Solder and connection 1-forms,    $\theta$ and $\omega$,   define a 1-form
\begin{equation}
\tau \equiv \theta \oplus  \omega   \, \colon ~ T\whM ~ \longrightarrow ~ \mR^q \oplus  \mathfrak{o}(q) ~ \cong  ~ \mR^{(q^2+q)/2}
\end{equation}
whose kernel is the distribution $T\whF$.

Recall that $\{E_1, \ldots , E_q\}$ is an orthonormal basis of $\mR^q$, and let $\{E_{ij} \mid 1 \leq i < j \leq q\}$ denote the corresponding  orthonormal  basis of $\mathfrak{o}(q)$.
Define orthonormal vector fields $\{\vec{Z}_1, \ldots, \vec{Z}_q\}$ on $\whM$ by specifying that at $\whx \in \whM$, 
\begin{equation}
\tau(\vec{Z}_i) = (E_i , 0) ~ , ~ \vec{Z}_i(\whx) \in  T_{\whx}\whF^{\perp}
\end{equation}
Similarly define vector fields $\{\vec{Z}_{ij} \mid  1 \leq i < j \leq q\}$ 
 \begin{equation}
\tau(\vec{Z}_{ij}) = (0,E_{ij}) ~, ~ \vec{Z}_{ij}(\whx) \in  T_{\whx}\whF^{\perp}
\end{equation}
The collection of  vector fields $\{\vec{Z}_i , \vec{Z}_{ij} \mid 1 \leq i < j \leq q\}$ span $ T_{\whx}\whF^{\perp}$ for each $\whx \in \whM$.

Recall that a    function $f \, \colon \whM \to \mR$ is $\whF$-basic if $f$ is constant on the leaves of $\whF$.   
Let $\cA = \cA(M,\F,g)$ denote the vector space consisting of    all linear combinations
 \begin{equation}
 \vec{Z} = \sum_{1 \leq i \leq q} ~ a^i \vec{Z}_i  ~ + ~ \sum_{1 \leq i < j \leq q} b^{ij} \vec{Z}_{ij}
  \end{equation}
where $\{a^1, \ldots , a^q\}$ and $\{b^{ij} \mid 1 \leq i < j \leq q\}$ are  $\whF$-basic, smooth functions on $\whM$.

Let $\cX(\whF)$ denote the smooth vector fields on   $\whM$ that are everywhere tangent to $\whF$.

One of the fundamental results of the Molino structure theory    is that  the flows of  the vector fields in $\cA$ are foliated:
\begin{prop} \label{prop-TP} Let $\vec{Z} \in \cA$ and 
  $\vec{X} \in \cX(\whF)$, then  $L_{\vec{Z}}(\vec{X}) \in \cX(\whF)$. Hence, for each $t \in \mR$, the flow of $\vec{Z}$, $\Phi^{\vec{Z}}_t \, \colon \whM \to \whM$,  maps leaves of  $\whF$ to leaves of $\whF$. $\hfill \eop$
\end{prop}
As the flows of vector fields in $\cX(\whF)$ preserve  the leaves of $\whF$, it follows that the group of foliated diffeomorphisms ${\bf Diff( \whM, \whF)}$  for $\whF$ acts transitively on $\whM$. That is, the foliated manifold $(\whM, \whF)$ is \emph{transversally complete} (TC), and the collection of  vector fields $\{\vec{Z}_i , \vec{Z}_{ij} \mid 1 \leq i < j \leq q\}$ define a \emph{transverse parallelism} (TP) for $\whF$.

\eject

Given  $\whx = (x,e) \in \whM$,   let $\whL_{\whx}$ denote  the leaf of $\whF$ containing $\whx$, 
and  $L_x$ be the leaf of $\F$ through $x$.  
Given a leafwise closed curved $\gamma \, \colon [0,1] \to L$ with $\gamma(0) = \gamma(1) = x$, 
we have a transverse holonomy map $h_x[\gamma]$ which depends only on the homotopy class of $\gamma$.  The differential  $dh_x[\gamma] \, \colon Q_x \to Q_x$   is an isometry, so induces a map $dh_x \, \colon \pi_1(L_x , x) \to \Isom(Q_x)$.   
Let $\cK_x \subset \pi_1(L_x , x)$ denote   the kernel of $dh_{x}$.

The framing $e \, \colon \mR^q \to Q_x$ induces an isomorphism $e^* \, \colon \Isom(Q_x) \cong \Oq$, and the composition $dh_{\whx} = e^* \circ dh_x \, \colon \pi_1(L_x , x) \to \Oq$ is the 
\emph{framed linear holonomy homomorphism} at $\whx$.

Given $\whx  = (x,e) \in \whM$, the leaf $\whL_{\whx}$ of $\whF$ is defined as an integral manifold of the flat connection on the $\Oq$-frame bundle over $L$, so that $y = (x,f) \in \whL_{\whx}$ means  that there is 
 $[\gamma] \in \pi_1(L_x , x)$ for which $f = dh_{x}[\gamma] (e)$.
The projection
 $\pi \, \colon \whL_{\whx} \to L_x$ is thus   the holonomy covering of $L_x$  associated to the kernel of $dh_x$. 

For simplicity of notation, let $\oL_{\whx}$ denote the closure $\owhL_{\whx}$  of a leaf $\whL_{\whx}$ of $\whF$. The distinction between $\oL_{\whx} \subset \whM$ and the leaf closure $\oL_x \subset M$ is indicated by the basepoint.

Recall that the foliation $\whF$ is invariant under the right action of $\Oq$ on $\whM$. For   $\whx \in \whM$,  define  two stabilizer subgroups of $\Oq$ associated to ${\whx}$:
\begin{eqnarray}
\cH_{\whx} & \equiv &  \{ A \in \Oq \mid \whL_{\whx}  \, A = \whL_{\whx} \}  \label{eq-iso1}\\
\bH_{\whx} & \equiv & \{A \in \Oq \mid \oL_{\whx}  \, A = \oL_{\whx}\} \label{eq-iso2}
\end{eqnarray}
Clearly, $\bH_{\whx}$ is the topological closure of $\cH_{\whx}$ in $\Oq$.

\begin{lemma} \label{lem-isotropy}
There is a natural identification of $\cH_{\whx}$ with the image of $dh_{\whx}$.     
\end{lemma}
\sop Let $\whx = (x,e)$ and  $[\gamma] \in \pi_1(L_x , x)$.  Set 
$$f \equiv dh_{x}[\gamma] (e) = dh_{x}[\gamma]  \circ e \in \Fr(Q_x)~ , ~   
A \equiv  e^{-1} \circ f =   e^{-1} \circ dh_{x}[\gamma]  \circ e \in \Oq$$
Then $\whx \, A = (x, e \, A) = (x, f)$. It follows that for $A = dh_{\whx}[\gamma]$, we have that $A \in \cH_{\whx}$.

Conversely, if $ \whL_{\whx}  \, A = \whL_{\whx}$ then $\whx \, A = (x, e \, A) = (x,f) \in \whL_{\whx}$ hence there exists $[\gamma] \in \pi_1(L_x , x)$ such that   $f  = dh_{x}[\gamma] (e)$. Thus, 
$A = e^{-1} \circ   dh_{x}[\gamma] \circ e$ is in the image of $dh_{\whx}$. $\hfill \eop$

\begin{cor} The following are equivalent:
\begin{enumerate}
\item $\cH_{\whx}$  is infinite
\item $\cH_{\whx} \subset \bH_{\whx}$ is a proper inclusion
\item  $\cK_x$ has infinite index in $\pi_1(L_x , x)$. $\hfill \eop$
\end{enumerate}
\end{cor}

  For $\epsilon > 0$, let  $\mD^q_{\epsilon} = \{ \vec{X} \in \mR^q \mid \| \vec{X} \| < \epsilon\}$. 
Let $Q^{\epsilon} \to M$ denote the unit disk subbundle, so that for each $\whx = (x,e) \in \whM$, 
the framing $e$ restricts to an isometry   $e \, \colon \mD^q_{\epsilon} \to Q^{\epsilon}_x$.

 The   lift  of $\nabla^L$  to $\whQ = \pi^*Q \to \whL_{\whx}$ is also flat, and  by the definition of $\whM$ the leaf $\whL_{\whx}$ has trivial holonomy.  Thus, combining the isometry  $e \, \colon \mD^q_{\epsilon} \to Q^{\epsilon}_x$ with    lifted Riemannian connection $\widehat{\nabla}^L$ on  $\whQ$   defines an isometric    product decomposition 
\begin{equation}
\xi_{\whx} ~ \colon ~ \whL_{\whx} \times \mD^q_{\epsilon} ~  \cong  ~ \widehat{Q}^{\epsilon} | \whL_{\whx}
\end{equation}

 Let $\exp \colon TM \to M \times M$ denote the geodesic exponential map, $\pi_2 \colon M \times M \to M$ the projection onto the second factor, 
 and  let $\exp_x = \pi_2 \circ \exp  \, \colon T_xM \to M$ be the exponential map based at $x \in M$.
 Choose $\epsilon > 0$ sufficiently small so that for all $x \in M$, the restriction 
$\exp_x \colon Q_x^{\epsilon} \to M$ is an embedding. 

\begin{prop}\label{prop-folprod}
The composition
\begin{equation}\label{eq-folprod}
\Xi_{\whx}   = \pi_2 \circ  \exp   \circ   d\pi   \circ \xi_{\whx} ~ \colon \whL_{\whx} \times \mD^q_{\epsilon} ~  \cong  ~ \widehat{Q}^{\epsilon}  ~  \longrightarrow  M
\end{equation}
is a foliated immersion. Given any $\why \in \whL_{\whx}$ and unit-vector $\vec{X} \in \mR^q$, the path
$$ \gamma_{(\whx, \vec{X})}(t) = \Xi_{\whx} (y, t  \, \vec{X}) ~, ~ -\epsilon < t < \epsilon$$
is a unit speed geodesic in $M$ which is orthogonal to  $\F$. $\hfill \eop$
\end{prop}
Note that Proposition~\ref{prop-folprod} does  not assert that the map 
 $\Xi_{\whx}$ is   an isometry, as the metric on $\mD^q_{\epsilon}$ is flat, while the   curvature tensor of $M$ transverse to $\F$ need not be zero. 

The fundamental  group $\pi_1(L_x , x)$ acts on the right on $\whL_{\whx}$ via   covering deck transformations, and   acts on $\mD^q_{\epsilon}$ via the holonomy representation $dh_x$.
 The product action on 
$\ds \whL_{\whx} \times \mD^q_{\epsilon}$ preserves the product structure, so we obtain a linear foliation 
$\F^{\omega}$ on the quotient 
\begin{equation} \label{eq-linear}
 Q^{\epsilon} | L_x  = \left( \whL_{\whx} \times \mD^q_{\epsilon}\right)/\pi_1(L_x , x) ~, ~ (\why \cdot \gamma, \vec{X})  \sim (\why , dh_x[\gamma] (\vec{X}))
\end{equation}
The map $\Xi_{\whx}$ is constant on the orbits of this action, so we obtain 
\begin{cor} \label{cor-locprod}
The induced map 
\begin{equation}\label{eq-locprod}
 \Xi_{\whx} ~ \colon ~ Q^{\epsilon} | L_x  \longrightarrow M
\end{equation}
 is a foliated immersion.  $ \hfill \eop$
\end{cor}

 Corollary~\ref{cor-locprod} implies  that a Riemannian foliation has a ``linear model'' in an open tubular neighborhood of a leaf. This is usually stated for the normal bundle to a compact leaf $L_x$ but is equally valid when formulated in terms of immersed submanifolds. The linear foliation  (\ref{eq-linear})  and map (\ref{eq-locprod}) yields a precise description of the leaves of $\F$   near to $L_x$.

 For $\whx =(x,e) \in \whM$, define the transverse disk to $\F$ as
 \begin{equation}\label{eq-transversal}
\iota_{\whx} ~ \colon ~ \mD^q_{\epsilon} \longrightarrow M ~, ~ \iota_{\whx}(\vec{X}) = \Xi_{\whx}(\whx , \vec{X})
\end{equation}
The image of $\iota_{\whx}$ will be denoted by $\cT^{\epsilon}_{x}$ which is a local transversal to $\F$ through $x$.  

Let $y \in \cT^{\epsilon}_{x}$ and $\vec{X} \in \mD^q_{\epsilon}$ so that $\iota_{\whx}(\vec{X}) = y$.
Then $\ds L_y$ is the image under   $\Xi_{\whx}$ of the leaf 
$( \whL_{\whx} \times \{\vec{X} \} )/\pi_1(L_x , x) \subset Q^{\epsilon} | L_x$.

\begin{prop} \label{prop-fixed}
The projection $\pi \, \colon Q^{\epsilon} | L_x \to L_x$ is a diffeomorphism when restricted to  $L_y$ for $y = \iota_{\whx}(\vec{X})$ if and only if 
$\vec{X}$ is fixed by all elements of $\cH_{\whx}$. 
\end{prop}
\sop
  $\cH_{\whx}$ is the image of the map $dh_x \, \colon  \pi_1(L_x , x) \to \Oq$ so that 
by (\ref{eq-locprod})   the covering map $\pi \, \colon L_y \to L_x$ has fibers isomorphic to the orbit   
$ \vec{X}  \cdot \cH_{\whx}$. 
$\hfill \eop$

\begin{cor} \label{cor-isotopy} Let $y \in \cT^{\epsilon}_{x}$. 
If $\pi \, \colon L_y \to L_x$ is a diffeomorphism, then there is a 1-parameter family of immersions, $I \, \colon L_x \times [0,1] \to M$, such that  $I_t \, \colon L_x \to M$ is a diffeomorphism onto a leaf of $\F$ for all $0 \leq t \leq 1$, $I_0 \colon L_x \to L_x \subset M$ is the inclusion of $L_x$,  and $I_1 \, \colon L_x \to L_y$.
\end{cor}
\sop For $z \in L_x$, define $I_t (z) =  \Xi_{\whx} (\{z\} \times \{t\vec{X}\})$, $0 \leq t \leq 1$. $\hfill \eop$

   \bigskip

   Finally, we recall the aspects of the     Molino structure theory for $\F$  which give a complete description of the closures of the leaves of $\F$, and of the foliation induced on them by $\F$.

\begin{thm}  \label{thm-molino1} Let  $\F$ bve a Riemannian foliation of a closed manifold $M$. \begin{enumerate}
\item For each $\whx \in \whM$, the leaf closure $\oL_{\whx}$ is a submanifold of $\whM$, and the set of all such leaf closures defines a foliation $\what{\E}$ of $\whM$ with all leaves compact.
\item For $\whx, \why \in \whM$ there exists a foliated diffeomorphism $\Phi_{\whx \why} \, \colon \whM \to \whM$ such that  $\Phi_{\whx \why}(\whx) = \why$,   hence  $\Phi_{\whx \why}(\whL_{\whx}) = \whL_{\why}$ and  $\Phi_{\whx \why}(\oL_{\whx}) = \oL_{\why}$ 
\item  There exists a closed manifold $\whW$ with a right $\Oq$-action, and an $\Oq$--equivariant fibration 
$\ds \whUp   \, \colon \whM \to \whW$ whose fibers of $\whUp$ are the  leaves of $\whE$. 
\item The metric on $T\whM$ defined above is projectable for the foliation $\what{\E}$, and for the induced metric on $T\whW$, the 
fibration $\ds \whUp   \, \colon \whM \to \whW$ is a Riemannian submersion. $\hfill \eop$
\end{enumerate}
\end{thm}

  Let $\ds W = M/\overline{\F}$ denote the Hausdorff space defined as the quotient of $M$ by the closures of the leaves of $\F$, and $\Up \, \colon M \to W$ the quotient map.  Then there is an $\Oq$-equivariant commutative diagram:
 
\begin{equation}
\begin{array}{rcccl}
&\Oq&{=}&\Oq&\\
 &\downarrow&&\downarrow& \\
&\whM&\stackrel{\whUp}{\longrightarrow}&\whW&\\
\pi&\downarrow&&\downarrow&\whpi\\
&M&\stackrel{\Up}{\longrightarrow}&W&
\end{array}
\end{equation}

Note that  given $\whx \in \whM$ and $\whw = \whUp(\whx)$, we have
$\ds \oL_{\whx} = \whUp^{-1}(\whw)$.

\begin{cor}\label{cor-leafclosure}
Let $\whx = (x,e)$ and $L_x$ be the leaf of $\F$ through $x$. Then 
\begin{equation}\label{eq-quotient}
\oL_x = \pi(\oL_{\whx}) = \oL_{\whx}/\bH_{\whx}
\end{equation}
The restriction   $\pi \, \colon \oL_{\whx} \to \oL_x$ is a principle $\bH_{\whx}$-fibration, and is a covering map if and only if  $\bH_{\whx}$ is a finite group. $\hfill \eop$
 \end{cor}

Let  $\cW$  be the horizontal distribution for the Riemannian submersion 
$\ds \whUp   \, \colon \whM \to \whW$. Then 
$\cW \subset   T\whM$ is the subbundle of vectors orthogonal to the fibers of $\ds \whUp   \, \colon \whM \to \whW$.  Note that $T\whF$  is contained in the kernel of $d\whUp$,   so that  $\cW \subset T\whF^{\perp}$ and each  leaf    $\whL_{\whx}$ is orthogonal  to the fibers of $\pi \, \colon \whM \to M$.

   The second part of the Molino structure theory concerns the geometry of  $\oL_{\whx}$ with the    foliation defined by the leaves of $\whF$.

\begin{thm}  \label{thm-molino2}   
   There exists a   connected, simply connected  Lie group $G$ with Lie algebra $\mathfrak{g}$ such that   the restricted foliation $\whF$ on $\oL_{\whx}$ is a Lie $G$-foliation with all leaves dense, defined by  a  $\mathfrak{g}$-valued connection 1-form 
$\ds \omega_{\mathfrak{g}}  \, \colon T \oL_{\why}  \longrightarrow  \mathfrak{g}$.
 
Moreover, given  $\whx \in \whM$ and a contractible  open neighborhood $\whV  \subset \whW$ of $\whw = \whpi(\whx)$ there exists  an  $\whF$-foliated diffeomorphism 
\begin{equation}\label{eq-product}
\Phi_{\whx} \, \colon \oL_{\whx} \times \whV  \longrightarrow \whU = \whpi^{-1}(\whV) \subset \whM
\end{equation}
Hence, $\whF | \whU$ is defined by a $\mathfrak{g}$-valued connection 1-form 
$\ds \omega_{\mathfrak{g}}^{\whU}  \, \colon T \whU  \longrightarrow  \mathfrak{g}$. $\hfill \eop$

\end{thm}

 \section{Equivariant foliated  transverse category}\label{sec-equivariant1}

In this section, we introduce the  \emph{$\Oq$-transverse category} of $\whF$ and show that this is equal to   the transverse category of $\F$.  The   proof uses the horizontal distribution  of a projectable metric on $M$ and Molino  theory, and is somewhat analogous to techniques used in the study of foliations with an Ehresmann connection \cite{BH1984a,BH1984b}.

Let  $\whU \subset \whM$ be an $\Oq$-invariant, $\whF$-saturated open subset. Let  
$\whH \, \colon \whU \times [0,1] \to \whM$ be an $\Oq$-equivariant,  $\whF$-foliated  homotopy. 
Then for $U = \pi(\whU)$,  $\whH$ descends to an $\F$-foliated homotopy   $H \, \colon U \times [0,1] \to M$.  The following result proves   the converse:

  \begin{prop} \label{prop-lift}
  Let  $H \, \colon U \times [0,1] \to M$ be an $\F$-foliated  homotopy. 
Then there exists an $\Oq$-equivariant, $\whF$-foliated  homotopy 
  \begin{equation}
\whH \, \colon \whU \times [0,1] \to \whM
\end{equation} 
such that $\pi \circ \whH = H \circ (\pi \times Id)$. That is, the following diagram commutes:
 
\begin{equation}
\begin{array}{rcccl}
&\whU \times [0,1]&\stackrel{\whH}{\longrightarrow}&\whM&\\
\pi \times Id&\downarrow&&\downarrow& \pi \\
&U \times [0,1] &\stackrel{H}{\longrightarrow}&M&
\end{array}
\end{equation}
  \end{prop}
  \sop
  Recall that   $\cH = \ker(\omega) \subset T\whM$ is the horizontal distribution  for the basic connection $\omega$.   For each $\why = (e,f) \in \whM$,   $d\pi \, \colon \cH_{\why} \to T_yM$ is an isomorphism.

     Given $\whx = (x,e) \in \whU$, let $\gamma_x(t) = H(x,t)$ denote  the path traced out by its homotopy.
Let  $\widehat{\gamma}_{\whx} \, \colon [0,1] \to \whM$ denote the horizontal lift of $\gamma_x(t)$ starting at $\whx$. That is,  $\widehat{\gamma}_{\whx}'(t) \in \cH$ and $\ds d\pi(\widehat{\gamma}_{\whx}'(t) ) = \gamma_x'(t)$      for all $0 \leq t \leq 1$. 
 Define $\whH(\whx, t) = \widehat{\gamma}_{\whx}(t)$.
  
  The map $\whH$ is smooth, as the lift of the path $\gamma_x$ to the solution curve 
  $\widehat{\gamma}_{\whx}(t)$ depends smoothly on the initial conditions $\whx$. 
  
  Given $A \in \Oq$ set $\why = (x,f) = (x,eA)$, then $R_A(\widehat{\gamma}_{\whx}(0))  = \why$, 
 \begin{eqnarray}
R_A(\widehat{\gamma}_{\whx}(t))' & = &  dR_a (\widehat{\gamma}_{\whx}'(t)) \in dR_A \cH = \cH \\
d\pi(dR_A (\widehat{\gamma}_{\whx}'(t))) & = &  (\pi(R_A(\widehat{\gamma}_{\whx})))'(t) = (\pi(\widehat{\gamma}_{\whx}))'(t) = \gamma_x'(t)
\end{eqnarray}
  so by uniqueness we have
  $R_A(\widehat{\gamma}_{\whx}(t)) = \widehat{\gamma}_{\why}(t)$. Thus, $\whH$ is $\Oq$-equivariant.

  It remains to show that $\whH$ is $\whF$-foliated.     
    
  Let $\cI^p = (-1,1)^p$, $\cI^q = (-1,1)^q$ and $\cI^m = (-1,1)^m$, where $m = p+q$. 
    
  Given $x \in U$, let 
  $\varphi \, \colon V \to \cI^m$ be a foliation chart such that $x \in V \subset U$. The  connected components of the leaves of  $\F | V$ are the plaques
  $$\cP_{\xi} = \varphi^{-1}(\cI^p \times \{\xi\}) ~, ~ \xi \in \cI^q$$
  Let $\varphi_{tr} \, \colon U \to \cI^q$ denote the projection onto the transverse coordinate $\xi$. 
  By assumption, the restriction of $g$ to  $Q | U$ projects to a Riemannian metric $g_U$ on $\cI^q$. 
 
Let  $\widehat{\cI^q} \to \cI^q$ denote   the     $\Oq$-bundle of orthogonal frames of $T\cI^q$, and 
  $\whpi \, \colon \widehat{\cI^q} \to \cI^q$ the projection.
Let
   $\omega_U \, \colon \widehat{\cI^q} \to \mathfrak{o}(q)$ be the Levi-Civita connection 1-form for $g_U$, with horizontal distribution $\cH_U \subset T \widehat{\cI^q}$.

The assumption that $g$ is projectable implies the restriction to $\whU =  \pi^{-1}(U)$  of the connection 1-form $\omega$ for $Q$  satisfies $\omega | \whU = \varphi_{tr}^* (\omega_U)$. In particular, this implies that the horizontal distribution $\cH$ contains the tangent vectors to the fibers of the projection map $\widehat{\varphi_{tr}} \, \colon  \whU \to  \widehat{\cI^q}$. Actually, this fact  is obvious as the fibers are exactly the plaques of $\whF | U$.
The following  technical result is used to show    that $\whH$ is a foliated map.    
\begin{lemma}   \label{lem-key}
 Assume  given smooth maps  $\widehat{G} \, \colon [0,\delta] \times [0,\epsilon] \to \whU$ and  $G \, \colon [0,\delta] \times [0,\epsilon] \to U$  such that $\pi \circ \widehat{G} =  G$. For 
   each $0 \leq s \leq \delta$, define the smooth    curves     $\sigma_s(t) = G(s,t)$  and     $\widehat{\sigma}_s(t) = \widehat{G}(s,t)$.
 Further assume that
 \begin{enumerate}
\item for each $0 \leq t \leq \epsilon$, the curve $\gamma_t$ defined by    $\gamma_t(s)  =G(s,t)$ is contained in a leaf of $\F | U$
\item the curve $\widehat{\gamma}_0$ defined by    $\widehat{\gamma}_0(s)  =\widehat{G}(s,0)$ is contained in a leaf of $\whF | \whU$
\item for all $(s,t) \in [0,\delta] \times [0,\epsilon] $, the tangent vector   $\widehat{\sigma}_s'(t) \in \cH | \whU$
\end{enumerate}
Then  the curve $\widehat{\gamma}_{t}$ defined by $\widehat{\gamma}_{t}(s)  =\widehat{G}(s,t)$ is contained in a leaf of $\whF | \whU$.
\end{lemma}
  \sop 
  Consider the diagram
  \begin{equation}
\begin{array}{crcclcccl}
 &  [0,\delta] \times [0,\epsilon] & \stackrel{\widehat{G}}{\longrightarrow}&\whU &  &\stackrel{\widehat{\varphi_{tr}}}{\longrightarrow}& & \widehat{\cI^q}&\\
& &   &\downarrow&\pi &  & &\downarrow& \whpi \\  
&(s,t) \in  [0,\delta] \times [0,\epsilon] & \stackrel{G}{\longrightarrow}&U   & & \stackrel{\varphi_{tr}}{\longrightarrow}& & \cI^q&
\end{array}
\end{equation}
  By (\ref{lem-key}.1) the composition $\tau_s(t) = \varphi_{tr} \circ G(s,t)$ is constant in $s$, and thus defines a smooth path $\tau  \, \colon [0,\epsilon] \to \cI^q, \tau(t) = \tau_s(t)$ for any choice of   $0 \leq s \leq \delta$.
  Let $\widehat{\tau} \, \colon [0,\epsilon] \to \widehat{\cI^q}$ be the lift of $\tau$ to a horizontal path with respect to $\omega_U$.  
  
  For   $0 \leq s \leq \delta$,  then     $\widehat{\sigma}_s'(t) \in \cH | \whU$  by (\ref{lem-key}.3), hence
  $d\widehat{\varphi_{tr}} ( \widehat{\sigma}'_s(t)) =  (\widehat{\varphi_{tr}} \circ  \widehat{\sigma}_s)'(t)$ is horizontal  for $\omega_U$.
  
The assumption $\pi \circ \widehat{G} =  G$ implies     $d\pi( \widehat{\sigma}'_s(t)) = \sigma_s'(t)$ hence 
  $d\whpi (\widehat{\varphi_{tr}} \circ  \widehat{\sigma}_s)'(t)) = \tau_s'(t) = \tau'(t)$.
 Thus,    the curve   $\widehat{\varphi_{tr}} \circ \widehat{\sigma}_s(t)$  is a horizontal lift of  $\tau(t)$, with initial point 
 $\widehat{\varphi_{tr}} \circ \widehat{\sigma}_s(0)$.
 
  The initial point  $\widehat{\varphi_{tr}} \circ \widehat{\sigma}_s(0)$ is independent of $s$ by (\ref{lem-key}.2), thus it  follows that the curve 
  $\widehat{\varphi_{tr}} \circ \widehat{\sigma}_s(t)$ is independent of $s$.

   That is, for all $0 \leq t \leq \epsilon$,  the curve  $s \mapsto \widehat{\sigma}_s(t)$ is contained in a fixed fiber of $\widehat{\varphi_{tr}}$, which is  a plaque of $\whF | \whU$, so contained in a leaf of $\whF$.
  $\hfill \eop$

  \bigskip

We now complete the proof of Proposition~\ref{prop-lift}. 
Let  $\whx \in \whU$ and $\why \in \whL_{\whx} \cap \whU$ be in the plaque containing $\whx$. We must show that $\whH_t(\whx)$ and $\whH_t(\why)$ lie in the same leaf for all $0 \leq t \leq 1$.

Choose a smooth path $\what{\gamma}_{\whx \why} \, \colon [0,1] \to  \whL_{\whx} \cap \whU$ with $\what{\gamma}_{\whx \why}(0) = \whx$ and $\what{\gamma}_{\whx \why}(1) = \why$. 
Define $x = \pi(\whx)$, $y= \pi(\why)$, and set $\gamma_{xy}(s) = \pi(\what{\gamma}_{\whx \why} (s))$  so that $\gamma_{xy}$ is a smooth path in $L_x \cap U$ from $x$ to $y$.

Compose these paths  with the given homotopy $H$ and its lift $\whH$, to obtain maps 
\begin{eqnarray*}
\whH(s,t) & = &  \whH(\what{\gamma}_{\whx \why} (s), t)\\
H(s,t) & = &  H(\gamma_{xy}(s), t)
\end{eqnarray*}
The maps $H_0$ and $\whH_0$ are   inclusions, so   $\whH(s,0) = \what{\gamma}_{\whx \why} (s)$ and  $H(s,0) = \gamma_{xy}(s)$ for $0 \leq s \leq 1$.
We are also given that $\pi \circ \whH = H$, so that $\pi \circ \whH(s,t) = H(s,t)$ for all $0 \leq s \leq 1$, $0 \leq t \leq 1$.

If the image of $H \, \colon [0,1] \times [0,1] \to M$ is contained in a foliation chart $\varphi \, \colon V \to \cI^p \times \cI^q$, then we can directly apply Proposition~\ref{prop-lift} to obtain the claim.  

For the general case, observe that there exists an integer $N > 0$ so that for 
$s_{\mu} = \mu/N$ and $t_{\nu} = \nu/N$ for each $0 \leq \mu, \nu < N$ there is a foliation chart 
$\varphi_{\mu \nu} \, \colon V_{\mu \nu} \to \cI^p \times \cI^q$
such that $H([s_{\mu}, s_{\mu+1}] \times [t_{\nu} , t_{\nu+1}]) \subset  V_{\mu \nu}$. 

Set $\delta = \epsilon = 1/N$, and define for $0 \leq s , t \leq 1/N$
$$
\widehat{G}_{\mu \nu}(s,t)   =   \whH(s_{\mu} + s, t_{\nu} + t) ~, ~  G_{\mu \nu}(s,t)   =    H(s_{\mu} + s, t_{\nu} + t) 
$$

 For $\nu = 0$ and each $0 \leq \mu < N$, the maps $\widehat{G}_{\mu 0}$, $G_{\mu 0}$ satisfy the hypotheses of 
Proposition~\ref{prop-lift}. The conclusion of  Proposition~\ref{prop-lift}   implies that this is again true for $\nu =1$, so that one can apply the Proposition repeated to obtain that the curve $s \to \whH(s,1)$ lies  in a leaf of $\whF$ as claimed.
$\hfill \eop$

  \bigskip

\begin{defn} Let  $\whU \subset \whM$ be an $\Oq$-invariant, $\whF$-saturated open subset. 
 We say that  $\whU$ is   \emph{$\Oq$-transversely  categorical} if there exists an $\Oq$-equivariant,  $\whF$-foliated  homotopy 
$\whH \, \colon \whU \times [0,1] \to \whM$ such that $\whH_0$ is the inclusion, 
and $\whH_1$ has image in the  orbit $\oL_{\whx} \cdot \Oq$ of the closure $\oL_{\whx}$ of a leaf $L_{\whx}$ of $\whF$. 
\end{defn}

The  \emph{$\Oq$-transverse category} of $\whF$, denoted by $\catO(\whM, \whF)$,  is the least number of $\Oq$-invariant, $\whF$-saturated open  sets required to cover $\whM$.  The results of this section imply:
  
\begin{cor}\label{cor-equality1}
Let $\F$ be   a Riemannian foliation of a compact manifold $M$, then 
\begin{equation}\label{eq-main2a}
\ecatt(M,\F) =  \catO(\whM,\whF) 
\end{equation}
\end{cor}

  \vfill
  \eject

\section{$\Oq$-equivariant  category}\label{sec-equivariant2}

In this section, we show that  the   $\Oq$-transverse category  $\catO(\whM, \whF)$ of $\F$     is equal to   the 
$\Oq$-equivariant  category  $\catO(\whW)$ of $\whW$, which will complete the proof of Theorem~\ref{thm-main3}. The   proof uses the horizontal distribution of the Riemannian submersion $\whUp \, \colon \whM \to \whW$ defined by  the projectable metric for $\whF$.

Let  $\whU \subset \whM$ be an $\Oq$-invariant, $\whF$-saturated open subset. Suppose that 
$\whH \, \colon \whU \times [0,1] \to \whM$ is an $\Oq$-equivariant,  $\whF$-foliated  homotopy. 
Then by Proposition~\ref{prop-induced} applied to $\whF$, we obtain an   $\Oq$-invariant  open set 
$\cU = \whUp(\whU) \subset \whW$ and  an $\Oq$-equivariant   homotopy  $\cH \, \colon \cU \times [0,1] \to \whW$ which is the quotient map of $\whH$.
Conversely, we have:

\begin{prop}\label{prop-equiv}
Let $\cU \subset \whW$ be     an $\Oq$-invariant  open   subset. Given   an $\Oq$-equivariant   homotopy  $\cH \, \colon \cU \times [0,1] \to \whW$, then for $\whU = \whUp^{-1}(\cU) \subset \whM$, 
 there is  an $\Oq$-equivariant, $\whF$-foliated  homotopy $\whH \, \colon \whU \times [0,1] \to \whM$  such that $\whUp   \circ \whH = \cH \circ (\whUp \times Id)$.
\end{prop}
\sop
For each $(\whw, t) \in \cU \times [0,1]$, 
set $\whw_t = \cH(\whw, t)$, and  let 
$$\cH'(\whw,t) = \frac{d}{dt}\cH(\whw,t) \in T_{\whw_t}\whW$$
 denote  the tangent vector along the time coordinate. 

For $A \in \Oq$, let $R_A \, \colon \whW \to \whW$ denote the right action  $R_A(\whw) = \whw  \, A$. 
The assumption that  $\cH$ is $\Oq$-invariant implies   
\begin{equation}
\cH'(\whw  \, A, t) =  \frac{d}{dt} \left( R_A \cH(\whw,t)\right) = dR_a \left( \cH'(\whw, t) \right)
\end{equation}
so that $\cH'(\whw,t)$ is a $\Oq$-invariant vector field.

Recall that  $\cW \subset T\whM$  is the subbundle  of vectors orthogonal to the fibers of $\whUp \, \colon \whM \to \whW$. Then $\cW$  is $\Oq$-invariant, as the Riemannian metric $\whg$ on $T\whM$ is $\Oq$-invariant, and $\whUp$ is $\Oq$-equivariant. 

For $\whx \in \whU$ with image $\whw = \whpi(\whx)$, define a smooth curve   \,
$\whx_t \, \colon [0,1] \to \whM$ by requiring that   
\begin{equation}\label{eq-lift2}
\whx_0 = \whx ~, ~ \whx'_t = \frac{d \whx_t}{dt} \in \cW ~ , ~ d\whUp (\whx'_t) = \cH'(\whw,t) =  \whw'_t 
\end{equation}
Thus, $\whx_t$ is an integral curve for the horizontal distribution  $\cW$. 
Define $\whH(\whx, t) = \whx_t$. 

It follows from (\ref{eq-lift2}) that     $\whUp(\whx_t) = \whw_t$, 
hence  $\whUp   \circ \whH = \cH \circ (\whUp \times Id)$.

The function $\whH$ is smooth as the integral curves $\whx_t$ depend smoothly on the initial condition $\whx_0 = \whx$. 
The function $\whH$ is   $\Oq$-equivariant, as given $A \in \Oq$  
$$   \frac{d}{dt} R_A(\whx_t)   = dR_A(\whx'_t) \in dR_A(\cW) = \cW$$
and 
$$ d\whUp\left( \frac{d}{dt} R_A(\whx_t) \right) 
= d\whUp \left(  dR_A(\whx'_t) \right) 
= dR_A \left( d\whUp (\whx'_t) \right) 
= \cH'(\whx  \, A, t) $$

It remains to show that $\whH_t \, \colon \whU \to \whM$ is  $\whF$-foliated. 
For each $\whx \in \whM$, the trace $t \mapsto \whx_t = \whH_t(x)$ is determined by the flow of the non-autonomous vector field $\whH' = \whH'(\whx, t)$ which is $\whUp$-related to the vector field $\cH' = \cH'(\whw , t)$.

Given a vector field  $\vec{X} \in \cX(\whF)$ tangent to $\whF$, let $\Phi_s^{\vec{X}} \, \colon \whM \to \whM$ denote its flow.  The metric $\whg$ on $T\whM$ is $\whF$-projectable, so its projection to $T\whF^{\perp}$ is invariant under the flow  $\Phi_s^{\vec{X}}$.

The flow $\Phi_s^{\vec{X}}$ preserves the leaf closures of $\whF$, so induces a quotient flow on $\whW$ which is constant. Set $\whx_{t,s} =  \Phi_s^{\vec{X}}(\whx_t)$, then $\whUp(\whx_{t,s})$ is constant as a function of $s$.

Apply  $\Phi_s^{\vec{X}}$  to the $\whUp$-related vector field $\whx_t' = \whH'(\whx, t)$ to obtain 
$$ \vec{Y}_{t,s} =  d\Phi_{s}^{\vec{X}}(\whH'(\whx, t))  =    d\Phi_{s}^{\vec{X}}(\whx'_t) $$
Note that $\frac{d}{ds}|_{s=0}\vec{Y}_{t,s} = L_{\vec{X}}H'(\whx,t)$. 
We claim that the Lie bracket 
$$(L_{\vec{X}}H')_{\whx_t} = [\vec{X}, \whH']_{\whx_t} \in T_{\whx_t}T\whF$$
which implies that the flow of $\whH'$ preserves the foliation $\whF$.

At each point $\whx_{t,s} \in \whM$ there is an orthogonal decomposition 
$$T_{\whx_{t,s}}\whM  ~ = ~ T_{\whx_{t,s}}^{\F}\whM ~ \oplus ~ T_{\whx_{t,s}}^{\whUp}\whM ~ \oplus ~  \cW_{\whx_{t,s}}$$
where 
  $T_{\whx_{t,s}}^{\F}\whM = T_{\whx_{t,s}}\whF$, and 
  $T_{\whx_{t,s}}^{\whUp}\whM =  T_{\whx_{t,s}}\oL_{\whx_{t,s}} \, \cap \, T_{\whx_{t,s}}\whF^{\perp}$ consists of the vectors orthogonal to
 $T_{\whx_{t,s}}\whF$ and tangent to the fibers of $\whUp$.  We decompose $\vec{Y}_{t,s}$ into its components, 
$$ \vec{Y}_{t,s} =   \vec{Y}_{t,s}^{\whF} + \vec{Y}_{t,s}^{\whUp} + \vec{Y}_{t,s}^{\cW} $$
  
The projection $d\whUp( \vec{Y}_{t,s})$ is constant as a function of $s$, so 
$$d\whUp(\frac{d}{ds}\vec{Y}_{t,s}^{\cW}) = \frac{d}{ds} d\whUp(\vec{Y}_{t,s}^{\cW}) =  \frac{d}{ds} d\whUp( \vec{Y}_{t,s}) = 0$$
and as $d\whUp \, \colon \cW_{\whx_{t,s}} \to T_{\whw}\whW$ is an isomorphism, we conclude that 
$\frac{d}{ds}|_{s=0}\vec{Y}_{t,s}^{\cW} = 0$.

The map $d\Phi_s^{\vec{X}}$ induces an isometry on $T\whF^{\perp}$, so the length of the vectors 
$\vec{Y}_{t,s}^{\whUp}$ is a constant function of $s$. At $s=0$ we have that 
$\vec{Y}_{t,0}^{\whUp} = \whH'(\whx,t)^{\whUp} =0 $ as $\whH'(\whx,t)  \in \cW_{\whx_t}$.
 Hence, $\vec{Y}_{t,s}^{\whUp} = 0$ for all $s$.  Thus, we have
 $$L_{\vec{X}}H'(\whx,t) =  \frac{d}{ds}|_{s=0}\vec{Y}_{t,s} =  
\frac{d}{ds}|_{s=0}   \vec{Y}_{t,s}^{\whF} + \frac{d}{ds}|_{s=0} \vec{Y}_{t,s}^{\whUp} + \frac{d}{ds}|_{s=0} \vec{Y}_{t,s}^{\cW} = \frac{d}{ds}|_{s=0}   \vec{Y}_{t,s}^{\whF} \in T_{\whx_t}\whF$$
This completes the proof of  Proposition~\ref{prop-equiv}.
 $\hfill \eop$

\bigskip
 
\begin{cor}\label{cor-equality2}
Let $\F$ be   a Riemannian foliation of a compact manifold $M$, then 
\begin{equation}\label{eq-main2b}
 \catO(\whM,\whF) =  \catO(\whW)
\end{equation}
\end{cor}

Note that the proof of Proposition~\ref{prop-equiv}  implies that 
$ \catO(\whM, \whF) =  \catO(\whM, \whE)$.  
This is equivalent to (\ref{eq-main2b}) as there is a  natural equivalence 
$ \catO(\whM, \whE) = \catO(\whW)$ 
because  $\whUp \, \colon \whM \to \whW$ is a fibration with compact fibers.

\eject

\section{$G$-equivariant category and orbit type}\label{sec-orbit}

  In this section, we recall some general properties  of a smooth action of a compact Lie group $G$ on a compact manifold and  their applications to   $G$-equivariant category. References for this material are 
  \cite{ALQ2001,Bartsch1993,Bredon1972,Colman2002b,DK2000,KM1989,Marzantowicz1989,tomDieck1987}. These results    will be applied to the   case of the $\Oq$-action on $\whW$ in the next section.

 Let $G$ be a compact Lie group,  and   $R \colon N \times G \to N$  a smooth right action
  on a closed manifold $N$ such that the quotient $N/G$ is connected.  For $A \in G$,  denote $u \, A = R(u,A)$, and let $u \, G = \{u \, A \mid A \in G\}$ denote the orbit of $u$. 
   Define the closed stabilizer  subgroup of the action of $G$ on $N$,   
\begin{equation}\label{eq-stabilizer}
G_0 = \{ A \in G \mid u \, A = u ~ {\rm for ~ all} ~ u \in N\}
\end{equation}
Note that $G_0$ is always normal. The action of $G$ is said to be  \emph{effective}  if $G_0$ is the trivial subgroup.

 Let $U \subset N$ be a $G$-invariant open set. A map $\cH \, \colon U \times [0,1] \to N$ is said to be a $G$-homotopy if $H$ is $G$-equivariant, and  $\cH_0 \, \colon U \to N$ is the inclusion. It is $G$-categorical if, in addition, $\cH_1 \, \colon U \to u  \, G$ has image in a single orbit, for some   $u \in N$.
A $G$-invariant subset $U \subset N$ is $G$-categorical if there exists a 
 $G$-categorical homotopy $\cH \, \colon U \times [0,1] \to N$. 
The $G$-category $\catG(N)$ of   $N$ is the least number of $G$-categorical open sets required to cover $N$.

 We recall some basic aspects of the geometry of a smooth $G$-action. We assume that $N$ has a $G$-invariant Riemannian metric $g_N$ on $TN$.

For $u \in N$,  $w \in u  \, G$, and $\epsilon > 0$, and  define the $\epsilon$-normal bundle  
\begin{eqnarray*}
\cN(w, \epsilon) & = &  \left\{   \vec{X} \in T_w(N)  \mid   \vec{X} \perp T_w(w  \, G) ~ , ~ \| \vec{X} \| < \epsilon \right\} \\
\cN(u  \, G, \epsilon) & = &  \left\{   \vec{X} \in T_w(N)  \mid  w \in u  \, G ~ , ~ \vec{X} \perp T_w(u  \, G) ~ , ~ \| \vec{X} \| < \epsilon \right\} 
\end{eqnarray*}

Denote the geodesic exponential map  by $\exp \, \colon TN \to  N$. For $w \in N$,  let $\exp_w \, \colon T_wN \to N$ be the exponential map based at $w$. 
Define the \emph{$\epsilon$-normal neighborhood} to $u  \, G$ by 
\begin{equation}\label{eq-normalnbhd}
\cU(u  \, G, \epsilon) = \left\{\exp_w(\vec{X})  \mid w \in u  \, G ~ , ~ \vec{X} \in \cN(w, \epsilon) \right\}
\end{equation}

For $w = u  \, A$, $A \in G$,    the derivative of the right action $R$ defines an isometric  linear representation of   the isotropy group, 
\begin{equation}
d_wR \, \colon G_w \to \Isom( T_w(u  \, G)^{\perp} )
\end{equation}

Note that  $d_u R_A \, \colon T_u(u  \, G)^{\perp} \to T_w(u  \, G)^{\perp}$  conjugates $d_w R$ to $d_u R$.

Let $\cN_0(w,\epsilon) \subset T_w(u  \, G)^{\perp}$ denote the  fixed vectors for the representation $d_wR$.

  \begin{thm}[Equivariant Tubular Neighborhood] \label{thm-GTN}
 For $u \in N$, there exists $\epsilon > 0$ so that the open   neighborhood $U = \cU(u  \, G, \epsilon)$  is a $G$-equivariant retract of $u  \, G$. 
 That is,     there exists a  $G$-equivariant homotopy $\cH \, \colon U \times [0,1] \to N$ such that  
  $\cH_t  | u \, G$ is the identity  for all  $0 \leq t \leq 1$,
 and  $\cH_1(U) = u \,  G$.
In particular, $U$ is   a $G$-categorical open neighborhood of $u  \, G$.
  \end{thm}
\sop For $\epsilon > 0$ sufficiently small,  
$\exp \, \colon \cN(u  \, G, \epsilon) \to  N$ is a diffeomorphism. 
The homotopy $\cH$ is then defined using the homothety 
\begin{equation}\label{eq-homothety}
 \cH(\exp_w(\vec{X}) , t) =    \exp_w(t  \, \vec{X})
\end{equation}
of  the normal geodesic map, where   $\vec{X} \in \cN(w,\epsilon)$.
$\hfill \eop$

  \medskip

The next  result is a generalization of Theorem~\ref{thm-GTN}, except that 
  there is no assertion that the action of $G$ on $U$ has a linear model.
  \begin{thm}[Equivariant Borsuk] \label{thm-GEB}
  Let $A \subset N$ be a closed, $G$-invariant subset.  Then  there exists a $G$-invariant  open neighborhood $A \subset U$ and a $G$-equivariant homotopy $\cH \, \colon U \times [0,1] \to N$ such that  
  $\cH_t | A$ is the identity for all $0 \leq t \leq 1$,
and  $\cH_1(U) = A$.
  \end{thm}

  If $H$ is a closed subgroup of $G$, we denote by $(H)$ the conjugacy class of $H$ in $G$.   While for an orbit $u  \, G$ the isotropy group
  $\ds G_{v} = \{ A \in G \mid v   A = v\}$ 
    depends on the choice of  $v \in u  \,  G$, the conjugacy class $(G_{v})$  does not, and is therefore an invariant of $u  \, G$. The conjugacy  class $(G_{u})$ is called the \emph{orbit type} of $u  \,   G$. 
  
We say that $u \, G$ is a \emph{principal orbit}  if $G_u = G_0$. 
An orbit $u \, G$ with  dimension less that that of a principal orbit is said to be a \emph{singular orbit}.  If $u \, G$ has the same dimension as a principal orbit, but the inclusion $G_0 \subset G_u$ is proper (and of finite index) then $u \, G$ is said to be an \emph{exceptional orbit}. 
  
One of the basic results for a smooth action  of a compact group on a compact manifold is that 
 it   has a finite set of orbit types, which follows easily from s is a consequence of Theorem~\ref{thm-GTN}. (See Proposition~1.2, Chapter IV  of \cite{Bredon1972}, or   Theorem 5.11 of \cite{tomDieck1987}). In the case of the action $R \, \colon N \times G \to N$, 
    there exists a finite collection of closed subgroups $\{G_0, \ldots G_k \}$ of $G$, where $G_0$ is   defined by (\ref{eq-stabilizer}),  such that 
    for all $u \in N$, there exists $\ell$ such that $(G_u) = (G_{\ell})$.

  There is a partial order on the set of orbit types of the $G$-space $N$:   for $u, v \in N$, 
\begin{equation}
(G_{u}) \leq (G_{v}) ~ {\rm if} ~ AG_{v}A^{-1} \subset G_{u} ~ {\rm for ~ some}~ A \in G
\end{equation}
We adopt the notation $[u] = (G_u)$, then    $[u] \leq [v]$ (resp. $[u] < [v]$)  means that the isotropy group of $v$ is conjugate to a (resp. proper) subgroup of the isotropy group of $u$. Thus, the orbit $v  \, G \cong G_v \backslash  G$
  is a fibration over the orbit  $u  \, G \cong G_u \backslash G$ and therefore should be considered ``larger''. Note that $[u] \leq (G_0)$ for all $ u \in N$.

Given a closed subgroup $H$ of $G$, define the \emph{$(H)$-orbit type subspaces}
\begin{eqnarray*}
N_{(H)} & = &  \{ u \in N \mid (G_u) = (H) \}\\
N_{\leq (H)} & = &  \{ u \in N \mid (G_u) \leq (H) \}
\end{eqnarray*}
where  $N_{(H)}$ is non-empty if and only if $H = G_{\ell}$ some $0 \leq \ell \leq k$;   define $N_{\ell} = N_{(G_{\ell})}$. 
 
Here are some  standard properties of the orbit type spaces; see \cite{Bredon1972,DK2000, tomDieck1987} for details.
\begin{thm} \label{thm-OTfibration} For $0 \leq \ell \leq k$, $N_{\ell}$ is a $G$-invariant submanifold, 
the quotient space $ N_{\ell} / G$ is a smooth manifold,  and 
the quotient map 
$$\pi_{\ell} \, \colon N_{\ell}  \longrightarrow N_{\ell} / G$$
is a fibration. 
Moreover, $\pi_{\ell}$ has the $G$-equivariant path lifting property: given  a smooth path $\sigma \colon [0,1] \to N_{\ell}/G$ with $\sigma(0) = v$,  there exists a $G$-equivariant smooth map $\Sigma \, \colon  v \, G \times  [0,1]  \to N_{\ell}$ such that $\Sigma_0 \, \colon v \cdot G \to N_{\ell}$ is the inclusion of the orbit $v \, G$, and $\pi_{\ell} \circ \Sigma(v \, A, t) = \sigma(t)$ for all $A \in G$. 
$\hfill \eop$
\end{thm}
 
 \begin{cor}\label{cor-GPLP}
 Given a smooth  path $\sigma \colon [0,1] \to N_{\ell}$ with $v= \sigma(0)$,  there exists a $G$-equivariant map $\Sigma \, \colon v \, G \times [0,1]  \to N_{\ell}$ such that $\Sigma_0$ is the inclusion, and 
 $\pi_{\ell} \circ \Sigma(v \, A, t) = \pi_{\ell} \circ \sigma(t)$. $\hfill \eop$ 
 \end{cor}
    
   There are various subtleties which arise in the study of the orbit type spaces $N_{(H)}$ for a smooth compact Lie group action. One is  that the quotient space $N_{\ell} / G$ need not be connected, and the connected components of $N_{\ell} / G$  need not all have the same dimensions. A second issue is that $N_{(H)}$ need not be closed, and the structure of the closure $\overline{N_{(H)}}$ is a fundamental aspect of the study of the action. We consider both of these points in the following. First recall (see, for example, Theorem~5.14 and Proposition~5.15 of \cite{tomDieck1987}):
   
  \begin{prop}\label{prop-POT} Recall that we assume $N/G$ is connected. Then
  the principal orbit space $N_0 = N_{(G_0)}$ is an open dense $G$-invariant submanifold of $N$ such that $N_0/G$ is connected. 
For $1 \leq \ell \leq k$,  the submanifold $N_{\ell} \subset N$ has codimension at least two. If there are no exceptional orbit types, then $N_0$ is connected.
  \end{prop}
      
The case where $N_0$ is not connected occurs when there exists $A \in G$ that acts as an involution of $N$ with a codimension-one fixed-point set.

      Next, we introduce the $\cZ$-stratification of $N$ associated to the orbit-type decomposition. 
 For $u \in N$,  let $\cZ_u \subset N$ denote the $G$-orbit of the connected component of $N_{(G_u)}$ containing $u$.  Note that for $u, v \in N$, either $\cZ_u \cap \cZ_v = \emptyset$ or $\cZ_u = \cZ_v$. Note that if  $G$ is connected, then $\cZ_u$ is also connected.
The inclusion   $u  \, G \subset \cZ_u$ can  be strict.  For example, this is always the case when $\cZ_u$ is not a closed   subset of $N$.

 \begin{prop}\label{prop-Gstrat}
 The collection  of sets $\cZ_u$ for $ u \in N$ form a finite stratification    of $N$.  
That is,  there exists a finite set of points $\{\eta_0 , \ldots , \eta_K\} \subset N$ such that for 
$\cZ_{i} = \cZ_{\eta_{i}}$
 \begin{equation}
N = \cZ_{0} \cup \cdots \cup \cZ_{K}
\end{equation}
As $N_0/G$ is connected, we can  require that  $\cZ_0 = N_0$. 
 \end{prop}
 \sop There exists a finite number of orbit types, and for each orbit type $(G_{\ell})$ the space   $N_{(G_{\ell})}$ is a finite union  of connected submanifolds.    $\hfill \eop$

The    \emph{$\cZ$-stratification} of $N$  is the collection of sets
\begin{equation}
\mathfrak{M}_G(N) =\left\{ \cZ_0, \cZ_1 , \ldots , \cZ_K \right\}
\end{equation}
  It satisfies the axioms of a Whitney stratification (see Chapter 2 of \cite{DK2000}.)  
  
 For the study of $G$-category, it is more natural to consider the $\cZ$-stratification of $N$ than the orbit-type stratification by the manifolds $N_{\ell}$. This is because a $G$-homotopy preserves   connected components of the orbit-type stratification, hence the $\cZ$-stratification captures more of the $G$-homotopical   information about the action.

\eject

A fundamental property of the $\cZ$-stratification is the incident relations between the closures $\overline{\cZ_{\ell}}$ of the strata. This motivates the following definition, which we will subsequently relate to the order-type relations between the strata.
Define the \emph{incidence partial order} on the collection of sets $\{\cZ_{1}, \ldots ,  \cZ_{K} \}$ by setting
\begin{equation}
\cZ_i  \lesssim \cZ_j ~ \Longleftrightarrow ~ \cZ_i \subset \overline{\cZ_j}
\end{equation}
Set $\cZ_i \approx \cZ_j$ if $\cZ_i  \lesssim \cZ_j$ and $\cZ_j  \lesssim \cZ_i $.
Note that $\cZ_i \approx \cZ_j$ implies that  $\cZ_i = \cZ_j$.

We require two fundamental technical lemmas, used to study the properties of the incidence partial order. The first      implies that the function $u \mapsto [u]$ is ``lower semicontinuous'' on $N$.
  \begin{lemma} \label{lem-lowersemicontinuous}
  For $u \in N$,  let $\cU(u  \, G, \epsilon)$ be an  $\epsilon$-normal neighborhood as in Theorem~\ref{thm-GTN}. 
   Then       $[u] \leq [v]$ for all $v \in \cU(u  \, G, \epsilon)$. Moreover,  $v \in  \cU(u  \, G, \epsilon)$ satisfies $[u] = [v]$   if and only if 
  $v = \exp_w(\vec{X})$ for some $w \in u  \, G$ and  $\vec{X} \in \cN_0(w,\epsilon)$.
    \end{lemma}
\sop For  $v \in \cU(u  \, G, \epsilon)$ there exists $w \in  u  \, A$  and   $\vec{X} \in \cN(w, \epsilon)$ such that $v = \exp_w(\vec{X})$. Then $B \in G_v$  if and only if $w = w \, B$ and $d_wR_B(\vec{X}) = \vec{X}$.   Let $A \in G$ such that $w = u \, A$, then $ B \in A^{-1} G_u A$ so $G_v \subset A^{-1} G_u A$ and hence $[u] \leq [v]$.
If $[u] = [v]$ then for all $ B \in A^{-1} G_u A$   we have that  $d_wR_B(\vec{X}) = \vec{X}$, hence 
$\vec{X} \in \cN_0(w,\epsilon)$.
Conversely, $\vec{X} \in \cN_0(w,\epsilon)$ implies $G_v = A^{-1} G_u A$ hence $[u] = [v]$. 
$\hfill \eop$

   The next   result implies that the function $t \mapsto [\cH(u,t)]$ is ``upper semicontinuous'' for $G$-equivariant homotopy.
     \begin{lemma} \label{lem-uppersemicontinuous} 
     Let $\cH \, \colon U \times [0,1] \to N$ be a $G$-homotopy.  For   $u \in U$, set $u_t = \cH(x,t)$.
  Then for $0 \leq t \leq 1$, we have
$\ds G_u \subset G_{u_t}$ and therefore  $[u_t] \leq [u]$.
\end{lemma}
\sop For $A \in G_u$ then 
$\ds u_t  \, A = \cH(u,t)  \, A = \cH(u  \, A,t) = \cH(u,t) = u_t$. $\hfill \eop$

\medskip

The following result establishes the relationships between the incidence partial order and the orbit-type partial order.
\begin{prop}\label{prop-incidence}
Let $v \in \cZ_i$ and $u \in \cZ_j$, for $i \ne j$ (and hence $\cZ_i \cap \cZ_j = \emptyset$.)

Suppose that $v \in  \overline{\cZ_j} - \cZ_j$. Then  $[v] < [u]$   and $\cZ_i \subset (\overline{\cZ_j} - \cZ_j)$, hence $\cZ_i \lesssim \cZ_j$.

In particular,   $\cZ_i \lesssim \cZ_j$  and $\cZ_i \not\approx \cZ_j$ implies that   $[v] < [u]$.
\end{prop}
\sop
Given    that $v \in \overline{\cZ_j} - \cZ_i$  there exists a sequence $\{u_{\ell} \mid \ell =1,2, \ldots\} \subset \cZ_j$ such that $\ds \lim_{\ell \to \infty} ~ u_{\ell} = v$. 
Let   $\epsilon > 0$ be such  that  $\cU(v  \, G, \epsilon)$ is a $G$-categorical   neighborhood of $v  \, G$. Then there exists $\ell$ such that $u_{\ell} \in \cU(v  \, G, \epsilon)$. 
Let  $w \in v  \, G$, $\vec{X} \in \cN(w, \epsilon)$ be such that $\exp_w(\vec{X}) = u_{\ell}$. 
Then by Lemma~\ref{lem-lowersemicontinuous}, 
$[v] < [u_{\ell}] = [u]$ unless $\vec{X} \in \cN_0(w, \epsilon)$. If $\vec{X} \in \cN_0(w, \epsilon)$, we have that $[v] = [u_{\ell}] = [u]$, and hence
$\ds \exp_w \left( \cN_0(w, \epsilon) \right) \subset \cZ_j$.
Thus, $v \in \cZ_j$ contrary to assumption, so we must have $[v] < [u]$.
 
 It remains to show that $\cZ_i \cap (\overline{\cZ_j} - \cZ_j) \ne \emptyset$ implies $\cZ_i \subset (\overline{\cZ_j} - \cZ_j)$.  This follows from an argument similar to the above. $\hfill \eop$.

\begin{cor}  \label{cor-closed}
  A stratum $\cZ_{j}$ is minimal for the incidence partial order if and only if $\cZ_{j}$ is a closed submanifold.
  \end{cor}
\sop If $\cZ_{j}$ is not closed, then there exists $v \in  \overline{\cZ_j} - \cZ_j$, and then  $v \in \cZ_i$ for some $i \ne j$. Then  $\cZ_i \lesssim \cZ_j$ by Proposition~\ref{prop-incidence}, so $\cZ_j$ is not minimal. The converse is obvious.
$\hfill \eop$

 The orbit-type function $u \to [u]$ is lower-semicontinuous by 
 Lemma~\ref{lem-lowersemicontinuous}, so it is natural to  also consider a notion of minimality based on continuity:
 
 \begin{defn} A stratum  $\cZ_{u}$ is said to be \emph{locally minimal} if there is an open $G$-invariant neighborhood $U$ of the closure $\overline{\cZ_{u}}$ such that $[u] \leq [v]$ for all $v \in U$.
  \end{defn}

\begin{prop}\label{prop-compact}
$\cZ_u$ is   locally minimal if and only if   $\cZ_u$ is closed. Hence,  $\cZ_u$ is locally minimal if and only if it  a least element for the incidence  partial order. 
\end{prop}
\sop
Assume that  $\cZ_u$ is closed, then by the Equivariant Borsuk Theorem~\ref{thm-GEB} there exists a $G$-invariant open neighborhood 
$U$ of $\cZ_u$ and a $G$-equivariant homotopy $\cH \, \colon U \times [0,1] \to M$ 
such that $\cH_1(U) = \cZ_u$. We claim that for all $v \in U$, $[u] \leq [v]$. 
Let $v \in U$, and set $v_t = \cH(v,t)$. Note that $v_1 \in \cZ_u$ so that  $[v_1] = [u]$. By Lemma~\ref{lem-uppersemicontinuous}, 
  $[v_t] \leq [v_0] = [v]$ for all $t$. In particular, $[u] = [v_1] \leq  [v_0] = [v]$, as was to be shown.

Conversely, assume that $\cZ_u$ is locally minimal, with $G$-invariant open neighborhood $U$ of $\overline{\cZ_u}$ as in the definition. If   there exists $v \in \overline{\cZ_u} - \cZ_u \subset U$, 
then $[v] < [u]$ by Proposition~\ref{prop-incidence}.  This contradicts the assumption that $v \in U$ satisfies $[v] \geq [u]$.  $\hfill \eop$

\begin{cor}\label{cor-locminset}
For each  $u \in N$, there exists a locally minimal stratum $\cZ_j \subset \overline{\cZ_u}$. $\hfill \eop$
\end{cor}

 \medskip

For the remainder of this section, we consider  the properties of $G$-homotopy with respect to the incidence and orbit-type partial orders, and give applications to $G$-category. The next result implies that $G$-homotopy preserves the $\cZ$-stratification.
 \begin{prop}\label{prop-stratinv}
 Let    $\cH \, \colon U \times [0,1] \to N$ be a $G$-homotopy.  
 Then for all $u  \in U$  and $0 \leq t \leq 1$,  
 \begin{equation}
\cH_t(U \cap \cZ_u) \subset  \overline{\cZ_u}
\end{equation}
 \end{prop}
\sop
    For   $u \in U$, set $u_t = \cH(x,t)$. Then  $u_0 = u \in U \cap N_{(G_{u})}$.
 Lemma~\ref{lem-uppersemicontinuous} implies that for all $0 \leq t \leq 1$, $[u_t] \leq [u]$, 
 hence $(G_{u_t}) \leq (G_u)$ so that $u_t \in N_{\leq (G_u)}$. 
  Define 
  $$s_0 = \sup \{s \mid u_t \in \overline{\cZ_u} ~ {\rm for ~ all} ~ 0 \leq t \leq s\}$$ 
 
 Note that $u_{s_0} \in \overline{\cZ_u}$. Suppose that $s_0 < 1$. Then    for all $ \delta > 0$ there exists $s_0 < t < s_0 + \delta$ such that $u_t \in N_{\leq (G_u)}$ but 
 $u_t \not\in \overline{\cZ_u}$. Set $w = u_{s_0}$.
 
  Lemma~\ref{lem-uppersemicontinuous} implies that 
  $G_u \subset G_{u_t}$ for all $0 \leq t \leq 1$, so in particular $G_u \subset G_w$.
Let $\cU(w \, G , \epsilon)$ be an $\epsilon$-normal neighborhood of $w \, G$, and 
$\exp_w \, \colon \cN(w,\epsilon) \to N$ the $\epsilon$-disk transverse to $w \, G$.

Let $\cN(w,\epsilon, G_u) \subset \cN(w, \epsilon)$ denote the vectors fixed by the subgroup 
$G_u$ under the    isotropy representation $d_w R \colon G_w \to \Isom( \cN(w,\epsilon))$. 
Then the submanifold $\exp_w (\cN(w,\epsilon, G_u))$ contains the intersection 
$\exp_w ( \cN(w,\epsilon)) \cap \cZ_u$ hence 
\begin{equation}
\exp_w (\cN(w,\epsilon, G_u)) =  \exp_w ( \cN(w,\epsilon)) \cap \overline{\cZ_u}
\end{equation}
 By assumption, there exists $s_0 < t < s_0 + \delta$ such that
 $u_t \in N_{\leq (G_u)}$ and there exists $\vec{X} \in \cN(w,\epsilon)$ but 
 $\vec{X} \not\in \cN(w,\epsilon, G_u) $ so that $v = \exp_w(\vec{X}) \in u_t \, G$.  
 Thus, $G_v \subset G_w$ as the elements of $G_v$ fix the vector $\vec{X}$, so we  have 
the proper inclusions  $G_u \subset G_v \subset G_w$. By hypothesis, $u_t \in \cZ_{\ell} \subset N_{\leq (G_u)}$    for some $\cZ_{\ell} \ne \cZ_u$ hence $G_u$ and $G_v$ are conjugate in $G$. This is impossible, as $G_u$ is a proper subgroup of $G_v$.
 $\hfill \eop$

 \begin{cor}\label{cor-locmin}
Let    $\cH \, \colon U \times [0,1] \to N$ be a $G$-homotopy.  Suppose that $\cZ_u$ is a locally minimal set for $u \in U$.  Then for all  $0 \leq t \leq 1$,  
 \begin{equation}
\cH_t(U \cap \cZ_u) =  \cZ_u
\end{equation}
  \end{cor}
\sop  By Corollary~\ref{cor-closed} the set $\cZ_u$ is closed, hence   $\cH_t(U \cap \cZ_u) \subset   \cZ_u$ by    Proposition~\ref{prop-stratinv}. As $\cZ_u$ is a closed submanifold, and $\cH_0$ is the identity,   the map  $\cH_t$ must be surjective for all $0 \leq t \leq 1$. $\hfill \eop$

 \bigskip
 Next, we consider the properties of the $G$-category and its relation to the $\cZ$-stratification.
\begin{prop}\label{prop-finiteGcat} The $G$-category 
$\Lambda = \catG(N)$ is finite. Moreover, there exists a $G$-categorical covering 
 $\{\cH_{\ell} \, \colon U_{\ell} \times [0,1] \to N \mid 1 \leq \ell \leq \Lambda\}$  and basepoints  $\{w_1, \ldots, w_{\Lambda}\} \subset N$   such that  each 
  $\cZ_{w_{\ell}}$ is locally minimal and $H_{\ell}(U_{\ell}) \subset w_{\ell}  \, G$.
  \end{prop}
\sop
By Lemma~\ref{lem-lowersemicontinuous}, every orbit has a $G$-categorical open neighborhood,  and $N$ compact implies there is a finite subcovering by $G$-categorical open sets, hence $\catG(N)$ is finite.

Let 
 $\{\cH'_{\ell} \, \colon U_{\ell} \times [0,1] \to N \mid 1 \leq \ell \leq \Lambda\}$  be a $G$-categorical open covering, with  $\cH'_{\ell}(U_{\ell}) \subset v_{\ell}  \, G$ for points 
   $\{v_1, \ldots, v_k\} \subset N$. 
   
   By Corollary~\ref{cor-locminset}, for each $1 \leq \ell \leq \Lambda$, we can  choose $w_{\ell} \in \overline{Z_{v_{\ell}}}$ such that $\cZ_{w_{\ell}}$ is a locally minimal set.      Let $\cU_{\ell} =  \cU(w_{\ell}  \, G, \epsilon)$ be an $\epsilon$-normal neighborhood 
   as in Theorem~\ref{thm-GTN}, with $G$-homotopy retract $\cH''' \, \colon \cU_{\ell} \times [0,1] \to N$.
(We choose   $\epsilon > 0$ sufficiently small so that it works for all $\ell$.) 
Then there exists $u_{\ell} \in \cU_{\ell} \cap Z_{v_{\ell}}$ in the same path-component as $v_{\ell}$.
    
    Let $\sigma_{\ell} \, \colon [0,1] \to Z_{v_{\ell}}$ be a smooth path such that $\sigma_{\ell}(0) = v_{\ell}$ and $\sigma_{\ell}(1) = u_{\ell}$.

By Corollary~\ref{cor-GPLP} each path $\sigma_{\ell}$ defines a smooth $G$-equivariant lifting 
$\Sigma_{\ell} \, \colon v_{\ell} \, G \times [0,1]  \to Z_{v_{\ell}}$ such that $\pi_{\ell} \circ \Sigma_{\ell}(v_{\ell} \, A , t) = \pi_{\ell} \circ \sigma_{\ell}(t)$. Thus, $\Sigma_{\ell}(v_{\ell} \, G , 1) \subset u_{\ell} \, G$.

Now, define $\cH_{\ell} \, \colon U_{\ell} \times [0,1] \to N$ as the concatenation:
$$
\cH_{\ell}(u,t)  = 
\left\{ 
\begin{array}{llccccccc}
        \cH'(u, 3t) &   {\rm for} & 0 & \leq & t & \leq & 1/3 , \\
        \Sigma_{\ell}(\cH'_{\ell}(u,1), 3t -1) &   {\rm for} & 1/3 & \leq & t & \leq & 2/3  \\
          \cH'''_{\ell}(\Sigma_{\ell}(\cH'_{\ell}(u,1),1), 3t -2) &    {\rm for} & 2/3 & \leq & t & \leq & 1   \\
\end{array}\right.
 $$
This yields a piece-wise smooth $G$-categorical homotopy $\cH$ as desired. 
By adjusting the time parameters, the map $\cH_{\ell}$ can be made    smooth.
$\hfill \eop$

\bigskip

One of the main problems for the study of $G$-category is to obtain upper and lower bounds for 
$\catG(N)$ in terms of the geometry and topology of the $G$-action.  The homotopy properties of   the $\cZ$-stratification  given above yields a geometric lower bound for $\catG(N)$.

Let  $\alpha_G(N)$ denote the number of locally minimum strata in   $\mathfrak{M}_G(N)$. 

 \begin{thm}\label{thm-est1}
$ \catG(N) \geq \alpha_G(N)$.
 \end{thm}
\sop
Let $\{\cH_{\ell} \, \colon U_{\ell} \times [0,1] \to N \mid 1 \leq \ell \leq k\}$ 
be a $G$-categorical covering of $N$, and let $\{w_1, \ldots, w_k\} \subset N$ be such that 
$H_{\ell}(U_{\ell}) \subset w_{\ell}  \, G$. 

Let $\cZ_i$ be a locally minimal set, and suppose that $u \in U_{\ell} \cap \cZ_i$. 
Set   $u_{\ell, t} = \cH_{\ell}(u,t)$.  Then   by Corollary~\ref{cor-locmin}, $u_{\ell, t} \in \cZ_i$ for all $0 \leq t \leq 1$, 
and hence   $u_{\ell, 1}  \, G = w_{\ell}  \, G \subset \cZ_i$. 

Suppose that $\cZ_i$ and $\cZ_i$ are disjoint locally minimal sets, such that 
    $u \in U_{\ell} \cap \cZ_i$ and $v \in U_{k} \cap \cZ_j$. 
By the above we have that 
$w_{\ell}  \, G \subset \cZ_i$ and $w_{k}  \, G \subset \cZ_j$. If $k = \ell$ then $\cZ_i \cap \cZ_j \ne \emptyset$, contrary to assumption. It follows that $k \ne \ell$.
Thus, for each locally minimal set $\cZ_i$ we can associate at least one $G$-categorical set $U_{\ell}$ so that distinct locally minimal sets yield distinct indices $\ell$. 
Hence $ \catG(N) \geq \alpha_G(N)$.
$\hfill \eop$

\bigskip

  Proposition~\ref{prop-stratinv} implies that the restriction of a $G$-categorical open set $U$ to a closed subset $\overline{\cZ_{\ell}}$ is again $G$-categorical. This remark yields another type of lower bound for $\cat_G(N)$:

 \begin{thm} \label{thm-est2} For each $\cZ_{\ell} \in \mathfrak{M}_G(N)$, 
$\ds   \catG(N) \geq  \catG(\overline{\cZ_{\ell}})$. 
 \end{thm}
 \sop A $G$-categorical covering of $N$ restricts to a $G$-categorical open covering of each closed subspace $\overline{\cZ_{\ell}}$.  
 $\hfill \eop$

In general, neither Theorems~\ref{thm-est1} or \ref{thm-est2} are optimal lower bounds, and best estimates are obtained by combining the ideas of each estimate for the particular  group action in question. 

 It is also possible to develop very sophisticated lower bound estimates for $\ds   \catG(N) $ in terms of the cohomology and homotopy theory of the compact group action \cite{CLOT2003}.

 There are also several types of upper bound estimates for $ \catG(N)$ as discussed for example  in the works of  
 Marzantowicz \cite{KM1989,Marzantowicz1989}, Bartsch \cite{Bartsch1993}, Ayala, Lasheras and Quintero  \cite{ALQ2001}, Colman \cite{Colman2002b} and the authors \cite{HT2006b}.
The simplest version is based on the dimension estimate, that for a connected   manifold
$X$ of dimension $\xi$, there is an upper bound $\cat(X) \leq \xi +1$. As each stratum  
$\cZ_{\ell} \in \mathfrak{M}_G(N) =\left\{ \cZ_0, \cZ_1 , \ldots , \cZ_K \right\}$ 
has quotient $\cZ_{\ell}/G$ which is a connected manifold, we can apply this estimate to each stratum to obtain

 \begin{equation}\label{eq-bestest}
 \alpha_G(N) ~  \leq ~   \catG(N) ~  \leq ~  \# \,  \mathfrak{M}_G(N)  ~  + ~  
\sum_{0 \leq \ell \leq K} ~ \dim \left( \cZ_{\ell}/G\right) 
\end{equation}

 \vfill
\eject

\section{Isotopy and $\cZ$-stratifications for $\F$}\label{sec-strata}

 In this section,  the results of the last section are    applied to the case  of a Riemannian foliation $\F$ and the associated space $\Oq$-manifold    $\whM$.   We first    introduce the \emph{isotopy stratification}  of $M$, corresponding to the structure of the leaves of $\F$. 
 Then the correspondence between the $\Oq$-orbits on $\whW$ and leaf closures for $\F$   is developed:  the main result is that  the $\cZ$-stratification of $\whW$  corresponds to  the  isotopy  stratification of $M$. This  yields an interpretation of the locally minimal sets $\cZ_u$ for the $\Oq$ -action on $\whE$  in terms of the intrinsic geometry of $\F$.   
References for this section are  the works of  Molino \cite{Molino1982,Molino1988,Molino1994},  Haefliger \cite{Haefliger1985,Haefliger1988} and Salem \cite{Salem1988}.

   \bigskip
   
 Let  $L$  be  a leaf of $\F$. An  \emph{$\F$-isotopy} of $L$ is a smooth map  $I \, \colon L \times [0,1] \to M$ such that $I_0 \colon L \to M$ is the inclusion of $L$, and for each $0 \leq t \leq 1$, $I_t \colon L \to M$ is a diffeomorphism onto a leaf $L_t$  of $\F$. We then  say that the   leaf $L_1$ is \emph{$\F$-isotopic} to $L_0 = L$, and write $L_0 \sim L_1$. For example, Proposition~\ref{prop-TP} implies that every leaf of $\whF$ is $\whF$-isotopic in $\whM$ to every other leaf of $\whF$. However, for       $\F$ this need not be true.
    
    Let $\cI_L$ denote the set of leaves of $\F$ which are $\F$-isotopic to $L$. For $L = L_x$ set $\cI_x = \cI_{L_x}$. The set of $\F$-isotopy classes of the leaves of $\F$ defines the 
    \emph{isotopy stratification}  of $M$.  

 \bigskip
 
Next, we use the Molino theory to define the holonomy stratification  for  $M$. 
  Recall that for  $\whx = (x,e) \in \whM$ with $\whw = \whUp(\whx)$, the fiber $\whUp^{-1}(w) = \oL_{\whx}$.  
  The projection $\pi \, \colon \whM \to M$ restricts to a covering map $\pi \, \colon \whL_{\whx} \to L_x$, and by Corollary~\ref{cor-leafclosure} we have that 
  $$\oL_x = \pi(\oL_{\whx}) \cong \oL_{\whx}/\bH_{\whx}$$
  where 
 $\ds  \bH_{\whx}   \equiv   \{A \in \Oq \mid \oL_{\whx}  \, A = \oL_{\whx}\} = \overline{\cH_{\whx}}$
by (\ref{eq-iso2}).

For $\whw = \whUp(\whx)  \in \whW$, the  $\Oq$-orbit $\whw  \, \Oq$ lifts to the $\Oq$-orbit of $\oL_{\whx}$, which again projects to the leaf closure $\oL_{x}$. Thus, each  orbit $\whw  \, \Oq$ corresponds to exactly one leaf $\oL_x$ of  $\E$.
Let $\Upsilon \, \colon M \to W$ denote the quotient map to the leaf space of $\E$. Then 
we have a commutative diagram of $\Oq$-equivariant maps:
\begin{equation}
\begin{array}{rcccl}
&\whM&\stackrel{\whUp}{\longrightarrow}&\whW&\\
\pi&\downarrow&&\downarrow&\whpi\\
&M&\stackrel{\Up}{\longrightarrow}&W&
\end{array}
\end{equation}

The set of leaves of $\F$ without holonomy form  an open  dense subset,  $M_0 \subset M$. 
 Define open dense subsets  $\whM_0 = \pi^{-1}(M_0)$ and $\whW_0 = \whUp(\whM_0) \subset   \whW$.  The leaves  in $M_0$ and $\whM_0$ are said to be    \emph{regular}, and the points of $\whW_0$ are   \emph{regular} orbits.

 Let $\whw = \whUp(\whx)$. Then  $\Oq_{\whw} =  \bH_{\whx}$. It is immediate from the definitions that  $\whUp(\why) \in \whW_{\Oq_{\whw}}$ if and only if  $\cH_{\why}$ is conjugate in $\Oq$ to a dense subgroup of $\bH_{\whx}$. 
 
 For a closed subgroup $H \subset \Oq$, $\whx = (x,e) \in \whM$ and $\whw = \whUp(\whx)$, set:
$$
 \begin{array}{ccccccc}
\whM_{(H)} & = &  \whUp^{-1}\left(\whW_{(H)}\right) & ; & M_{(H)} & = & \pi\left( \whM_{(H)}  \right)\\
\whM_{\leq (H)} & = &  \whUp^{-1}\left(\whW_{\leq (H)}\right) & ; & M_{\leq (H)} & = & \pi\left( \whM_{\leq (H)}  \right)\\
\whZ_{\whx} & = &  \whUp^{-1}\left(\cZ_{\whw}\right) & ; & Z_x & = & \pi \left( \whZ_{\whx} \right)
\end{array}
$$
 Note that as $\whUp$ has connected fibers, $\whZ_{\whx} $ can also be described as the $\Oq$-orbit  of the connected component of  $\whM_{(\Oq_{\whx})}$ containing $\whx$, and $Z_x$ is   the connected component of $M_{(\Oq_{x})}$ containing $x$.
 
Let $\ds  \mathfrak{M}_{\F}(M) =\left\{ \cZ_0, \cZ_1 , \ldots , \cZ_K \right\}$
  be the  $\cZ$-stratification of $\whM$ for the action of $\Oq$.
Then there  exists a finite set of points $\{z_1, \ldots , z_K\} \subset M$ such that 
$Z_{x_i} = \pi(\whUp^{-1}(\cZ_{i}))$. Set $Z_i = Z_{x_i}$.

  A stratum $Z_{i}$ is said to be  
\emph{locally minimal}  if the    set $\cZ_{i}$ is locally minimal.

   \begin{prop} \label{prop-holostrat} 
   If   $\why   \in \whZ_{\whx}$, then   $\cH_{\why}$ is conjugate in $\Oq$ to  $\cH_{\whx}$.
   \end{prop}
   \sop 
   Given  $\why = (y,f)   \in \whZ_{\whx}$,   either  there exists a continuous path $\sigma \, \colon [0,1] \to 
 \whM_{(\Oq_{\whx})}$ such that $\sigma(0) = \whx$ and $\sigma(1) = \why$, or  there exists $A \in \Oq$ such this holds for $(y,f  \, A)$.  So without loss of generality we can assume that $\why$ is in the same path component as $\whx$. 
 
 Let $\epsilon > 0$ be as in 
      Corollary~\ref{cor-locprod}. For each $\why_t = \sigma(t)$ there exists a linear model for $\F$ along the leaf $L_{y_t}$ through the point $y_t = \pi(\why_t)$. The image $\pi(\sigma[0,1])$ is compact, so is covered by a finite collection of such linear models. 
          It thus suffices to consider the case where  $x$ and $y$ are such that 
there is a foliated immersion as given in        Corollary~\ref{cor-locprod} 
 $$ 
 \Xi_{\whx} ~ \colon ~ Q^{\epsilon} | L_x = 
 \left( \whL_{\whx} \times \mD^q_{\epsilon}\right)/\pi_1(L_x , x) \longrightarrow M
$$
and $y \in \cT^{\epsilon}_{x} \cap  Z_{x}$. Let $\vec{X} \in \mD^q_{\epsilon}$ so that $\iota_{\whx}(\vec{X}) = \Xi_{\whx}(\what{\xi}) = y$, where $\what{\xi} = (\{ \whx \} \times \{ \vec{X} \}) \in Q^{\epsilon} | L_x$. 
The map $\Xi_{\whx}$ defines an orthonormal framing of $Q_y$ which we denote by $f$, so that   $\why = (y, f) \in \whM$.

By Lemma~\ref{lem-isotropy}, 
the holonomy group  of the foliation $\F^{\omega}$ on $Q^{\epsilon} | L_x$ at the point 
 $\what{\xi}$  is given by
\begin{equation}\label{eq-restricted}
\cH^{\omega}_{\what{\xi}} = \left\{ A \in \cH_{\whx} \mid \vec{X}  \, A = \vec{X} \right\} \subset \cH_{\whx} \subset \Oq
\end{equation}
Thus  $\cH^{\omega}_{\what{\xi}}$ is conjugate to 
 $\cH_{\why}$ in $\Oq$.

Let $\cV_{\whx} \subset  \mD^q_{\epsilon}$  be the linear subspace of vectors fixed by the subgroup $\cH_{\whx}$.
Note that  $\cV_{\whx}$ is also the set of vectors fixed by the closure 
$\bH_{\whx} = \overline{\cH_{\whx}}$.

Similarly, define $\cV_{\what{\xi}} \subset  \mD^q_{\epsilon}$ as the  linear subspace of vectors fixed by 
$\cH^{\omega}_{\what{\xi}}$.
Again,  $\cV_{\what{\xi}}$  is also the set of vectors fixed by the closure 
$\ds \overline{\cH^{\omega}_{\what{\xi}}} \subset \Oq$ which is conjugate to $\bH_{\why}$.

The key to the proof of   Proposition~\ref{prop-holostrat} is the following result for linear actions.
\begin{lemma}\label{lem-linearstruct} ~
 $\ds \cH^{\omega}_{\what{\xi}}  = \cH_{\whx}$ ~ $\Longleftrightarrow$ ~ $\cV_{\whx} = \cV_{\what{\xi}}$ ~ 
 $\Longleftrightarrow$ ~  $ (\bH_{\whz}) = (\bH_{\whx})$.
\end{lemma}
\sop
First note that we always have that 
$\ds \cH^{\omega}_{\what{\xi}}  \subset \cH_{\whx}$, hence 
$\ds  \overline{\cH^{\omega}_{\what{\xi}}} \subset    \bH_{\whx}$ 
and so $\ds \cV_{\whx} \subset \cV_{\what{\xi}}$.
If we show that $\ds \cH^{\omega}_{\what{\xi}}  = \cH_{\whx}$ then the other equalities follow immediately.

If $\cV_{\whx} = \cV_{\what{\xi}}$, then by (\ref{eq-restricted}) the vector $\vec{X} \in \cV_{\whx}$ hence 
$\ds \cH^{\omega}_{\what{\xi}}  = \cH_{\whx}$.

Since $\cV_{\whx}$ and $\cV_{\what{\xi}}$ are relatively closed subspaces,    $\cV_{\whx} \subset \cV_{\what{\xi}}$ is proper  implies the inclusion 
$\ds \bH_{\whx} \subset  \overline{\cH^{\omega}_{\what{\xi}}}$ is proper. Hence, 
$\ds \bH_{\whx} =  \overline{\cH^{\omega}_{\what{\xi}}}$ implies that $\cV_{\whx} = \cV_{\what{\xi}}$.
$\hfill \eop$

 \medskip
 We conclude the proof of Proposition~\ref{prop-holostrat}. 
 The assumption that $\why   \in \whZ_{\whx}$ implies  $ (\bH_{\whz}) = (\bH_{\whx})$. By Lemma~\ref{lem-linearstruct} we have that $\ds \cH^{\omega}_{\what{\xi}}  = \cH_{\whx}$. As remarked above, 
 $\ds \cH^{\omega}_{\what{\xi}}$ is conjugate to $\cH_{\whx}$.
   $\hfill \eop$

 \bigskip

 Proposition~\ref{prop-holostrat} and its proof yield the following characterization of the strata. 
     \begin{prop}\label{prop-holochar} For all $x \in M$, 
      $Z_{x} = \cI_x$. Hence, the isotopy stratification of $M$ is finite.
   \end{prop}
     \sop We first show that  $Z_{x} \subset  \cI_x$.
   Let $\whx = (x,e) \in \whM$ and let $y \in  Z_{x}$. 
   As in the proof of  Proposition~\ref{prop-holostrat}, it will suffice to show this for the case where 
  there is a foliated immersion
   $$ 
 \Xi_{\whx} ~ \colon ~ Q^{\epsilon} | L_x = 
 \left( \whL_{\whx} \times \mD^q_{\epsilon}\right)/\pi_1(L_x , x) \longrightarrow M
$$
and $y \in \cT^{\epsilon}_{x} \cap  Z_{x}$. Let $\vec{X} \in \mD^q_{\epsilon}$ so that $\iota_{\whx}(\vec{X}) = \Xi_{\whx}(\what{\xi}) = y$, where $\what{\xi} = (\{ \whx \} \times \{ \vec{X} \}) \in Q^{\epsilon} | L_x$. 
The map $\Xi_{\whx}$ defines an orthonormal framing of $Q_y$ which we denote by $f$, so that   $\why = (y, f) \in \whM$. 

Then by Proposition~\ref{prop-holostrat} we have that 
 $\ds \cH^{\omega}_{\what{\xi}}  = \cH_{\whx}$, 
 hence  $\vec{X} \in   \cV_{\whx}$ by the proof of  Lemma~\ref{lem-linearstruct}. 
 
Define an isotopy $I \, \colon L_x \times  [0,1] \to M$ where for $\whz \in L_{\whx}$ with $\pi(\whz) = z$, 
\begin{equation}
I(z, t) =  \Xi_{\whx} \left( \{\whz \} \times \{t\vec{X}\} \right) 
\end{equation}
where the map is well-defined as   $\vec{X} \in \cV_{\whx}$ and the action of $\pi_1(L_x , x)$ on $\cV_{\whx}$ is trivial. For the same reason, for all $0 \leq t \leq 1$,  the restriction 
$I_t \, \colon L_x \to M$ is a diffeomorphism onto the leaf through $\iota_{\whx}(t \vec{X})$.
As $\iota_{\whx}(\vec{X}) = y$ and $\iota_{\whx}(0) = x$, this proves that $y \in \cI_x$.

   Conversely, to show that  $\cI_x \subset  Z_{x}$, it suffices to consider an isotopy 
      $I \, \colon L_x \times  [0,1] \to M$ with $I_0(x) = x$ and $I_1(x) = y$, whose image lies inside a     foliated immersion
   $$ 
 \Xi_{\whx} ~ \colon ~ Q^{\epsilon} | L_x = 
 \left( \whL_{\whx} \times \mD^q_{\epsilon}\right)/\pi_1(L_x , x) \longrightarrow M
$$
The image $L_t = I_t(L_x)$ is isotopic to $L_x$ so   the  covering maps $L_t \to L_x$ are diffeomorphisms, thus the  holonomy groups $\cH_{\whx_t}$ have constant conjugacy classes. Hence     the conjugacy class of the  closure  $(\bH_{\whx_t})$ is constant, so $L_t \subset Z_x$ for all $0 \leq t \leq 1$. Thus, $y \in Z_x$. 
   $\hfill \eop$

     \begin{cor} \label{cor-locminisotopy}
     For $\whx = (x,e) \in \whM$, set $\whw = \whUp(\whx)$.  Then  $\cZ_{\whw}$ is locally minimal if and only if $\cI_x$ is a closed submanifold of $M$. $\hfill \eop$
   \end{cor}
 
    \begin{cor}\label{cor-duality}
 For each  isotopy class $\cI_x$, there is a well-defined conjugacy class $(\cH_x)$ defined as the conjugacy class of the holonomy group $\cH_{\whx}$ for any $\whx = (x,e) \in \whM$. $\hfill \eop$
 \end{cor}

  Corollary~\ref{cor-duality} highlights one of the unique aspects of the study of Riemannian foliations, the duality between intrinsic and extrinsic geometry of leaves.  That is, for $L_y$ near to $L_x$, we have that $L_y \sim L_x$ exactly when $L_y$ is diffeomorphic to $L_x$ via the orthogonal projection along leaves, and also that $L_y \sim L_x$ exactly when the  holonomy group $\cH_{\why}$ of $L_y$ is naturally conjugate to that of $L_x$.

\section{Foliations with finite category}\label{sec-finite}

Let $\F$ be   
    a Riemannian foliation $\F$ of a compact manifold $M$, then 
  the   $\Oq$-equivariant LS category  $\catO(\whW)$ is   finite, and by     
    Corollaries~\ref{cor-equality1} and \ref{cor-equality2}  we have the equalities        
$$
\ecatt(M,\F) =  \catO(\whM,\whF) =  \catO(\whW)
$$

In this section, we prove the promised   geometric criteria for when    $\catt(M,\F)$ is   finite, and show  that   $\ecatt(M,\F) = \catt(M,\F)$ when this criteria is satisfied. 
  Theorem~\ref{thm-main1}   follows from Propositions~\ref{prop-compactmin} and \ref{prop-equal}.
 As a further consequence, we  obtain   estimates for   $\catt(M, \F)$ which extend those proven  in \cite{CH2004} for compact Hausdorff foliations.

\begin{prop} \label{prop-compactmin}
Let $Z_x$ be a locally minimal set for $\F$. Suppose there is given  an open saturated set $U$ such that such that $U \cap Z_x \ne \emptyset$, and a foliated homotopy 
$H\, \colon U \times [0,1] \to M$ such that $H_1 \, \colon U \to M$ has image in a single leaf of $\F$.  Then every leaf of $\F$ in $Z_x$ is compact.
\end{prop}
\sop
Let $y \in U \cap Z_x$.  By Proposition~\ref{prop-holochar} every leaf of $\F$ in $Z_x$ is isotopic, so it suffices to show that there is some $z \in Z_x$ with $L_z$ compact. 

Let $z = H_1(y)$, so the image of $H_1 \, \colon U \to M$ is contained in the leaf $L_z$ of $\F$. 
As noted in section~2, the closure $\oL_y$ is a compact minimal set for $\F$.   By Lemma~\ref{lem-closures},   we have that  $\oL_{y} \subset U$, hence by   
Theorem~1.3 of \cite{Hurder2006a}, 
the image $H_1(\oL_y)$ is contained in a compact leaf  of $\F$, thus $L_z$ is compact. The proof   now follows from the following:
\begin{lemma}\label{lem-stratainv}
Let $Z_x$ be a locally minimal set for $\F$ and $U$ an open saturated set such that such that $U \cap Z_x \ne \emptyset$. If 
$H\, \colon U \times [0,1] \to M$ is  a foliated homotopy, then $H_t(Z_x \cap U) \subset Z_x$ for all $0 \leq t \leq 1$.
\end{lemma}
\sop
Let $\whw \in \whW$ be such that $Z_x = \pi(\whUp^{-1}(\cZ_{\whw}))$, and note that  $\cZ_{\whw}$ is a locally minimal set. 
 Let $\whU = \pi^{-1}(U) \subset \whM$, and   set $\cU = \whUp(\whU) \subset \whW$.

 Given $y \in  U \cap Z_x $, choose a point $\what{\xi} \in \cU \cap \cZ_{\whw}$  so that $\whUp(\pi^{-1}(y)) = \what{\xi}  \cdot \Oq$.

By Proposition~\ref{prop-lift},   the homotopy  $H$ determines an $\Oq$-equivariant homotopy $\cH \, \colon \cU \times [0,1] \to \whW$. 
  By Corollary~\ref{cor-locmin},    the trace $\what{\xi}_t = \cH(\what{\xi}, t) \in \cZ_{\whw}$ for all $0 \leq t \leq 1$,  and thus  $\pi(\whUp^{-1}(\what{\xi}_t)) \subset Z_x$.
 In particular, $L_z = \oL_z = \pi(\whUp^{-1}(\what{\xi}_1)) \subset Z_x$.
$\hfill \eop$

\begin{cor}\label{cor-fincpt}
Let $\F$ be     a Riemannian foliation $\F$ of a compact manifold $M$ with $\catt(M,\F)$   finite.  Then every local minimal set $Z_i$ for $\F$ consists of compact leaves.
\end{cor}
\sop
Choose    $x \in Z_i$. Then there a categorical open set $U$ with $x \in U$,  so by Proposition~\ref{prop-compactmin} every leaf of $Z_i$ is compact.
$\hfill \eop$ 
 
 \medskip

 \begin{prop}  \label{prop-equal} 
 Let $\F$ be     a Riemannian foliation $\F$ of a compact manifold $M$ such that  every local minimal set $Z_i$ for $\F$ consists of compact leaves. 
Then  $\catt(M, \F) = \ecatt(M,\F)$. 
\end{prop}
\sop
Let $\{\cH_{\ell} \, \colon \cU_{\ell} \times [0,1] \to \whW \mid 1 \leq \ell \leq k\}$ 
be an $\Oq$-categorical covering of $\whW$, for $k = \catO(\whW)$. Let $\{\whw_1, \ldots, \whw_k\} \subset \whW$ be such that 
$\cH_{\ell}(\cU_{\ell}) \subset \whw_{\ell}  \, \Oq$. 

By Proposition~\ref{prop-finiteGcat} we can assume that $\cZ_{\whw_{\ell}}$ is locally minimal for each $1 \leq \ell \leq k$. Thus, the sets $Z_{x_{\ell}}$ are locally minimal for $\F$. By assumption, all leaves in 
$Z_{x_{\ell}}$ are compact.

 For each $1 \leq \ell \leq k$,  let $\whx_{\ell}  \in \whUp^{-1}(\whw_{\ell})$  and set $x_{\ell} = \pi(\whx_{\ell})$ and $L_{\ell} = L_{x_{\ell}} \subset M$. Then $L_{x_{\ell}} \subset Z_{x_{\ell}}$ hence is a compact leaf.
 
 Let   $H_{\ell} \, \colon U_{\ell} \times [0,1] \to M$ be the $\F$-foliated homotopy corresponding to $\cH_{\ell}$ for  $U_{\ell} = \pi(\whUp^{-1}(\cU_{\ell}))$. 
 
 We have that  $H_1 \, \colon U_{\ell} \to \oL_{\ell}$, and as each $L_{\ell}$ is a compact leaf,   we are done.
$\hfill \eop$

\bigskip

More than just characterizing when $\catt(M,\F)$ is finite, the arguments of the previous sections and the above  yield an estimate for the transverse LS category. 
Proposition~\ref{prop-holochar} identifies the isotropy stratification of $\F$ with the $\cZ$-stratification of the $\Oq$-action on $\whW$. Let $\{I_i, \ldots , I_K\}$ be an enumeration of the isotropy strata for $\F$, and assume that $I_{\ell}$ is a closed submanifold exactly when $1 \leq \ell \leq k$,  where $k \leq K$. 
Set $\alpha(M,\F) = k$. Then by Lemma~\ref{lem-stratainv},  Theorem~\ref{thm-est2} and equation (\ref{eq-bestest}) we obtain

\begin{thm}\label{thm-bestest2}
Let $\F$ be a Riemannian foliation of a compact manifold $M$. Then
 \begin{eqnarray}
 \alpha(M,\F)  ~  \leq ~   \ecatt(M,\F)  &  \leq &   \sum_{1 \leq \ell \leq K} ~  \ecatt(I_{\ell} , \F | I_{\ell})  \label{eq-strataest}\\
&  \leq &~  K + \sum_{1 \leq \ell \leq K} ~ \dim \left( I_{\ell}/\overline{\F_{\ell}} \right) \label{eq-dimest}
\end{eqnarray}
Moreover, if every locally minimal set $Z_i$ consists of compact leaves, then $\catt(M,\F) = \ecatt(M,\F)$ so the estimates (\ref{eq-strataest}) and (\ref{eq-dimest}) also hold   for  the transverse LS category.
\end{thm}

Theorem~\ref{thm-bestest2} is a complete generalization of the estimate given by  Theorem~6.1 in \cite{CH2004} for the tranverse LS category of compact Hausdorff foliations.

\vfill
 \eject
 
\section{Critical points}\label{sec-critical}
 
 One of the applications of the LS-category invariant for a compact manifold $N$ is to give a lower bound estimate on the number of critical points for a $C^1$-function $f \colon N \to \mR$
(see \cite{CLOT2003,Fadell1985,FH1987,James1995,LS1934,Palais1966b,PT1988}.) 
 When there is a compact Lie group $G$ acting on $N$ and the function $f$ is invariant for the $G$-action, then there is an induced map on the quotient space, 
 $\overline{f} \colon N/G \to \mR$, 
  and one can attempt to estimate the number of critical points of $f$ using $\overline{f}$. Unfortunately, the quotient space $N/G$ need not be a manifold, so the classical theory does not apply. In addition, examples show that the category of $N/G$ can be much smaller than the category $\cat_G(N)$.
 
The solution is to consider   the set of critical points for $f$ as a $G$-space, and then the  $G$-category $\cat_G(N)$ provides a lower bound for the number of critical $G$-orbits \cite{Bartsch1993,Fadell1985,FH1987,Marzantowicz1989}.
 
Let  $f \colon M \to$ be      a $C^1$-function on a compact manifold $M$ with a Riemannian foliation $\F$. The function $f$ is said to be  $\F$-basic if it is constant along leaves of $\F$. Since $f$ is continuous, it is constant on the  leaves of the SRF $\E$ of $M$ defined by the closures of the leaves of $\F$, thus  induces a continuous map on the quotient space $\phi = \overline{f} \colon W = M/{\E} \to \mR$.   The differential $df \, \colon TM \to \mR$ is a basic 1-form, so if $L$ is a critical, then $\oL$ will consist of critical leaves also. Colman studied in section 5 of \cite{CM2001} the relation between the transverse category $\catt(M,\F)$ of $\F$ and the number of critical leaves of $\E$  for $f$, in the case where  all leaves of $\F$ are compact, so it is a compact Hausdorff foliation. 
 
An alternate approach is to first lift $f$ to a smooth function 
$\what{f} = f \circ \pi \, \colon \whM \to \mR$ which is basic for $\whF$.  
A leaf $\whL_{\whx}$ of $\whF$ will be critical for $\what{f}$ if and only if $L_x$ is critical for $f$, 
and so its closure $\oL_{\whx}$ is a critical submanifold for $d\what{f}$.  
Moreover,  the function $\what{f}$ is    $\Oq$-invariant, 
so we can estimate the number of critical leaf closures for $f$ in terms of the $\Oq$-invariant critical submanifolds of $\whM$.  The smooth map $\what{f}$ descends to a smooth $\Oq$-invariant map 
$\what{\phi} \colon \whM \to \mR$ and the $\Oq$-invariant critical sets for $\what{\phi}$  correspond exactly to the critical leaf closures of $f$. Now, the quotient space $\whW$ is a manifold, so we can apply the usual results of equivariant LS-category theory to estimate the number of critical $\Oq$-orbits for $\what{\phi}$. This is a great advantage, as the technical estimates required for the theory can be done in the context of a compact group action, instead of the case of a foliation with non-compact leaves. We recall the main result (see \cite{Fadell1985,FH1987,Marzantowicz1989,PT1988}):

  \begin{thm}   Let $G$ be a compact Lie group,  and   $R \colon N \times G \to N$  a smooth right action
  on a closed manifold $N$.  If $\what{\phi} \, \colon N \to \mR$ is a $C^1$-function which is $G$-invariant, 
  then  
  $\cat_G(N)$ is a lower bound on the number of critical orbits for $\what{\phi}$.
 \end{thm}
  
  Theorem~\ref{thm-main3} and the above discussion then yields the claim of Theorem~\ref{thm-LSR} as a consequence:
  \begin{cor}
   If $f \, \colon M \to \mR$ be a $C^1$-map which is constant along the leaves of $\F$, 
 then $\ecatt(M,\F)$ is a lower bound for the number of critical leaves of $\E$.
  \end{cor}

\bigskip

\section{Examples}\label{sec-examples}

In this section, we present a collection of examples to illustrate the ideas of the paper.  
There are three general methods   for constructing a Riemannian foliation on a compact manifold: isometric Lie group actions; the group suspension construction applied to an isometric action of a  finitely generated group; and the various blow-up constructions for singular Riemannian foliations 
\cite{Alex2004, AT2005, Molino1988, Molino1994}. Note that for open manifolds, there is an important  fourth method, which realizes a Riemannian pseudogroup as a Riemannian  foliation of  an open manifold \cite{Haefliger1971,Hurder1981}. This will  not be discussed, as little is known  of the transverse LS-category for open manifolds.

\begin{ex}\label{ex-compact}
Compact Hausdorff foliations and finite group actions
\end{ex}

Epstein \cite{Epstein1976} and Millett    \cite{Millett1974} proved that a foliation $\F$ with all leaves compact and whose leaf space     $M/\F$ is   Hausdorff   with the quotient topology, admits a projectable Riemannian metric, hence is Riemannian.  For each $x \in M$, the holonomy $\cH_x$ of the leaf $L_x$ is always a finite group, hence all leaves of the lifted foliation $\whF$ of the orthonormal frame bundle $\whM$ are also compact. The lifted foliation $\whF$ is thus defined by the fibration $\whUp \, \colon \whM \to \whW$.  

The quotient space $W = M/\F$ is a Satake manifold \cite{Satake1956}, or generalized orbifold, as every point is modeled either on $\mR^q$, or by a quotient of $\mR^q$ by a finite isometry group given by the holonomy of the leaf fiber. The quotient map $\whpi \, \colon \whW \to W$  is an ``$\Oq$-desingularization''  of $W$, where the regular orbits of $\Oq$ correspond to the leaves of $\F$ without holonomy. 
A leaf $L_x$ with holonomy for $\F$ corresponds to an orbit  $\whw \cdot \Oq$ on $\whW$ with  isotropy group  $\bH_{\whw}$ that strictly contains the stabilizer group of the action.  The foliation $\whF$ has no singular leaves.

The {\it exceptional set} $E_{\F}$ of a compact foliation $\F$ is   the union of all  leaves with holonomy, and thus corresponds to the union of all strata except for $\cZ_0$, hence 
$\ds  E_{\F} = \cZ_1 \cup \cdots \cup \cZ_K$.
 The   set of leaves without holonomy,   $G_{\F} =  M - E_{\F}$, is  called  the {\it good set}
and corresponds to the set of regular orbits for the $\Oq$-action on $\whW$.
The exceptional set $E_{\F}$ admits the \emph{Epstein filtration} by the holonomy groups of its leaves, and Proposition~\ref{prop-holostrat} shows that  the connected components of the strata in the Epstein filtration correspond to the isotropy stratification of $M$.

Theorem~5.3 of \cite{CH2004} proves that the Epstein filtration is invariant under $\F$-foliated homotopy. This result is a direct  consequence of Proposition~\ref{prop-stratinv} of this paper.
Moreover, Theorem~\ref{thm-bestest2}  above  implies the estimate of Theorem~6.1 in \cite{CH2004}.  
In fact, the genesis of this current work was to find a new approach, which would work for all  Riemannian foliations,  of the results for compact Hausdorff foliations in \cite{CH2004}.  

\medskip

{\bf Example \ref{ex-compact}.1}: ~
The papers \cite{Colman2002b,Colman2006,CH2004} contain constructions and calculations of the transverse LS category for compact Hausdorff foliations, and the reader is referred to those papers for details.  The standard  method of construction  is to start with a fibration $\wtpi \, \colon \wtM \to \wtW$ of  compact manifolds, and assume that $\wtpi$ is equivariant with respect to a finite group $\Gamma$ which acts freely on $\wtM$. Then the foliation of $\wtM$ by the fibers of $\wtpi$ descends to a compact Hausdorff foliation of $M$, whose leaf space $W = \wtW/\Gamma$ is thus a good orbifold.  Colman has also given examples of compact Hausdorff foliations whose leaf space $M/\F$ is a bad orbifold.

 \medskip 
 
\begin{ex}\label{ex-isometric}
Isometric flows
\end{ex}

Let $(M,g)$ be a compact, connected Riemannian manifold, and  $\phi \, \colon M \times \mR^p \to M$ a non-singular  isometric action of $\mR^p$. (The case $p=1$ corresponds to an  isometric flow on $M$.)  The orbits of $\phi$ define a Riemannian foliation $\F_{\phi}$  for which the metric  $g$ is   projectable.  The geometric and topological properties of this class of Riemannian foliations had been studied by many authors \cite{Carriere1984,Epstein1984,Ghys1984,HS1990,HS1991}. The leaves of the foliation $\whF$ are   given by a free $\mR^p$ isometric action on $\whM$. The closure of the image of $\mR^p$ in $\Isom(\whM)$ is a torus $\mT^k$ for some $k > p$, and  the foliation $\whE$ is defined by a free   isometric action of $\mT^k$ on $\whM$. This class of examples reveals many of the properties of the transverse category theory for Riemannian foliations, and we give  
 several examples to illustrate various phenomena.
 
 \medskip
 
 {\bf Example \ref{ex-isometric}.1}: ~ The canonical example of an ``irrational flow on the torus'' is formulated generally as follows. 
 Let $M = \mT^n$ be the $n$-torus, considered as the quotient $\mT^n = \mR^n/\mZ^n$ by the integer lattice.  Choose  $1 \leq p < n$ and a  real  matrix $A^{p \times n}$ such that $A A^T$ is invertible.  Then $A$ defines an injective map $\bA \colon \mR^p \to \mR^n$,  and thus yields an isometric affine action $\phi_A \, \colon \mR^p  \times \mT^n \to \mT^n$. 
 The orbits of $\phi_A$ are the affine planar leaves of $\F_A$.  All leaves of $\F_A$ are regular, as there is no holonomy. 

If all entries of $A$ are rational numbers, then the leaves of $\F_A$ are compact tori; otherwise, the leaves have closures which are embedded tori $\xi_A \, \colon \mT^k \subset \mT^n$ for some $p < k \leq n$. Thus,  $\catt(\mT^n, \F_A) = \infty$ unless $A$ is a rational matrix.. On the other hand, the essential transverse category is equal to the category of the foliation $\E_A$ obtained from the closures of the leaves.
 The leaf space $\mT^n/\E_A \cong \mT^n/\xi_A(\mT^k)$ is   a torus of dimension $n-k$, hence  
 $\ecatt(\mT^n , \F_A) = \cat(\mT^{n-k}) = n-k + 1$.
 
 \medskip

{\bf Example \ref{ex-isometric}.2}: ~ The previous examples of isometric flows can be embedded  into compact space forms. We begin with the simplest examples of this.

Let $\vec{\alpha} = (\alpha_0, \alpha_1, \cdots , \alpha_n) \in \mR^{n+1}$ and let $M$ be the unit $2n+1$ sphere
$$M = \mS^{2n+1} = \{ x = [z_0, z_1, \ldots, z_n] \mid |z_0|^2 + \cdots |z_n|^2 = 1 \}$$ 
Define an isometric $\mR$-action on $\mS^{2n+1}$ by 
$$
\phi_t([z_0, \ldots, z_n]) = [e^{2 \pi \alpha_0 t  \sqrt{-1}} z_0, \ldots , e^{2 \pi \alpha_n t  \sqrt{-1}} z_n]
$$
The orbits of $\phi_t$ define the leaves of the foliation $\F_{\vec{\alpha}}$.
 With the assumption that the numbers $\{1,\alpha_0 , \ldots , \alpha_n\}$ are linearly independent over $\mQ$, then the leaves of $\overline{\F_{\vec{\alpha}}}$ are defined by the action of the compact abelian group $\mT^{n+1} = \mS^1 \times \cdots \mS^1$ acting   diagonally.

 There are precisely $n+1$ locally minimal $\F_{\vec{\alpha}}$-isotopy classes,  corresponding to the orbits of the points $\vec{e}_i = [0, \ldots, 1, \ldots , 0]$, which are isolated circles. 
 Thus, $\catt(\mS^{2n+1}, \F_{\vec{\alpha}})  \geq n+1$.
 It is also easily seen that for each point $\vec{e}_i$ there is a flow-equivariant retraction of the open set 
$\ds \cU_i = \{ [z_0, \ldots, z_n] \in \mS^{2n+1} \mid z_i \ne0 \}$
to the orbit of $\vec{e}_i$. Hence  $\catt(\mS^{2n+1}, \F_{\vec{\alpha}})  = n+1$.
 
 It is possible to construct a wide variety of variations on this example, based on the general setup where $G$ is a connected, compact Lie group and $K \subset G$ is a closed subgroup, 
and we set $M = G/K$. Let $n$ denote the $\mR$-rank of $G$, so there is a locally free action
$\Phi_t \, \colon \mR^n \times G \to G$ whose orbit through the identity $e \in G$ is a maximal torus $\mT^n \subset G$. 
 
  Choose  $1 \leq p \leq n$ and a  real  matrix $A^{p \times n}$ such that $A A^T$ is invertible.  Then $A$ defines an injective map $\bA \colon \mR^p \to \mR^n$,  and thus yields an isometric affine action $\phi_A \, \colon \mR^p  \times \mT^n \to \mT^n$. 
  
For $\vec{v} \in \mR^p$ and $x \in G/K$,  define the action 
 $\phi_A(x , \vec{v}) = \Phi_t(x, \phi_A(\vec{v}))$ so that we obtain a $p$-dimensional foliation $\F_A$ on the homogenous space $G/K$.

The closure of the orbit of $\phi_A$ through the identity is a compact $k$-torus $\mT^k_A \subset \mT^n$  for some $p < k \leq n$. The foliation $\E_A$  --  defined by the closures of the leaves of $\F_A$ -- is given by the orbits of the closed subgroup $\mT^k_A$ on $G/K$.
 
 Note that in this generality, there is no assurance that $\F_A$ has any compact leaves, hence generically one has  $\catt(G/K, \F_A) = \infty$. On the other hand,  $\ecatt(G/K, \F_A)$ equals  the equivariant category of $G/K$ for the left action of the compact Lie subgroup $\mT^k_A  \subset G$. 
The calculation of $\ecatt(G/K, \F_A)$  then follows by methods of the theory of compact Lie groups,  and  has applications to  the  residue theory for the secondary classes of $\F_A$  \cite{HT2006b,LP1976a,LP1976b, Yamato1979, Yamato1981}.

 \medskip

\begin{ex}\label{ex-products}
Products
\end{ex}

 There is a simple remark, that if $M_1, \F_1)$ and $(M_2 , \F_2)$ are Riemannian foliations of compact manifolds with leaf dimensions $p_1$ and $p_2$ respectively, then the product manifold $M = M_1 \times M_2$ has a Riemannian foliation $\F = \F_1 \times \F_2$ whose leaves have dimension $p = p_1 + p_2$. As an example, suppose that $M_2 = \mT^{n_2}$ has a linear foliation as in Example \ref{ex-isometric}.1 with all leaves dense. Then for any Riemannian foliation $(M_1 , \F_1)$ the product foliation $\F$ has no compact leaves, hence 
 $\catt(M, \F) = \infty$. On the other hand, $\ecatt(M, \F) = \ecatt(M_1 , \F_1)$.

\medskip

\begin{ex}\label{ex-CLG}
Suspension and compact Lie group actions
\end{ex}
Given an action of a finitely generated group  on a compact $q$-dimensional manifold, 
$\alpha \, \colon \G \times N \to N$, 
the suspension construction yields a foliation $\F_{\alpha}$ of codimension $q$ whose transverse holonomy group is globally defined by the action (for example, see \cite{CN1985,CandelConlon2000}.) When the given action is isometric, then this yields a Riemannian foliation, and provides a large class of examples. We use this construction to realize the orbit structure of every compact Lie group action as the transverse geometry of some Riemannian foliation.

\medskip

{\bf Example \ref{ex-CLG}.1}: ~ 
 Let $\{\gamma_1, \ldots , \gamma_d\}$ be a set of generators for $\G$. A left isometric action of $\G$ on $N$ is equivalent to a   representation  $\alpha \, \colon \Gamma \to \Isom(N)$,  where each $\alpha(\gamma_i)$ acts isometrically on the left on $N$.
 
Let $B$ be a compact connected Riemannian manifold with basepoint $b_0 \in B$ 
such that the fundamental group $\Lambda = \pi_1(B, b_0)$ admits a surjection 
$\beta \, \colon \Lambda \to \G$. For example, one can let $B = \Sigma_d$ be a closed Riemann surface with genus $d$. 
There is a surjection $\Lambda \to \mF^d = \mZ * \cdots * \mZ$ 
onto the non-abelian free group on $d$ generators. 
The choice of the generators for $\G$ defines a surjection $\mF^d \to \G$, 
and the composition $\beta \, \colon \Lambda \to \G$ is then a surjection.

The universal cover $\wtB \to B$ has a right isometric action by $\Lambda$,  acting via deck transformations. 

Consider the product manifold $\wtB \times N$ with the product foliation $\widetilde{\F}$ whose leaves are the ``horizontal slices'' $\wtB \times \{x\}$ for $x \in N$. The Riemannian metrics on $\wtB$ and $N$ define the product metric   on $\wtB \times N$.

Define an action of $\Lambda$ on $\wtB \times N$ by specifying, for $\lambda \in \Lambda$,
$\ds  (b,x) \cdot \lambda = (b \cdot \lambda, \alpha \circ \beta (\lambda^{-1}) \cdot x)$. 

Both the product foliation and the product Riemannian metric on $\wtB \times N$ are invariant under this action, so the product foliation descends to a Riemannian  foliation $\F_{\alpha}$ of 
$\ds M = \wtB \times_{\Lambda} N$.  Note that there is an embedding $\iota_0 \colon N \to M$ given by $\iota(x) = [b_0 , x]$ where $[b,x]$ represents the equivalence class of the pair $(b, x)$ in $M$.  Let $N_0 \subset M$ denote the image of this map.

The group of isometries $\Isom(N)$ is compact, so the   closure 
$G = \overline{\alpha(\G)} \subset \Isom(N)$ is a compact Lie subgroup.
Let $\E_{\alpha}$ denote the singular Riemannian foliation of $M$ by the closures of the leaves of $\F_{\alpha}$. Then the leaves of $\E_{\alpha} \cap N_0$ are precisely given by the orbits of $G$.

\medskip

{\bf Example \ref{ex-CLG}.2}: ~ 

Let $G$ be a compact connected Lie group, and $\varphi \, \colon G \times N \to N$ a smooth isometric action on a compact Riemannian manifold $N$ of dimension $q$. The orbits of $\varphi$ define a singular Riemannian foliation $\E_{\varphi}$ of $N$ \cite{HS1990, HS1991,Molino1988,Molino1994}.  Ken Richardson   \cite{Richardson2001} showed that there always exists   a Riemannian foliation $\F_{\varphi}$ of a compact manifold $M$ such that 
the singular Riemannian foliation defined by the closures of the leaves of $\F_{\varphi}$ is transversally equivalent to $\E_{\varphi}$. We recall this argument.

Let $\G \subset G$ be a finitely generated dense subgroup; such   always exists  by a clever argument of Richardson   \cite{Richardson2001}. The restriction of $\varphi$ to the subgroup defines a representation 
$\alpha \, \colon  \G \to \Isom(N)$.

Use the  suspension construction as in Example \ref{ex-CLG}.2, to obtain a Riemannian foliation 
  $(M, \F_{\alpha})$ such that the image $\alpha(\G) \subset \Isom(N)$ has closure precisely the compact Lie subgroup $G$.    The orbit type stratification of the $G$-action on $N$ equals the stratification of the transversal $N_0$ induced by the closures of the leaves of $\F_{\varphi}$.

There is a particular case of this construction which yields some very interesting examples. Let $G = \bS\bU(n)$ be the group of $n\times n$ special unitary matrices.  Let  $N = \bS\bU(n)$ be the group itself, and let the action $\varphi$ be the adjoint, so that 
$\varphi(A) \, \colon \bS\bU(n) \to \bS\bU(n)$ is given by $\varphi(A)(B) = A^{-1} B A$. Let $\F_{\alpha}$ denote   the resulting suspension  foliation. Then $\ecatt(M, \F_{\alpha})  = n$, based on the calculations of \cite{HT2006b}. Note that the codimension of $\F_{\alpha}$ is the dimension of $\bS\bU(n)$ so that $q = n^2 -1$.

The suspension of the adjoint action of $\bS\bU(n)$ on itself yields a foliation  with no exceptional orbits. This is not the case with the groups $\bS\bO(n)$. In fact, for this case, there are isolated exceptional  orbits, so by Colman's results, the transverse category of the suspended foliation will be infinite. However, the essential transverse category  will be finite, and its calculation is a very important problem, as it yields estimates for the category of the groups $\bS\bO(n)$ themselves \cite{HT2006b}.

 \vfill



\begin{thebibliography}{10}


\bibitem{Alex2004}
{M.M.~Alexandrino},
\newblock {\it Singular {R}iemannian foliations with sections},
\newblock {\bf Illinois J. Math.}, 48:1163--1182, 2004.

\bibitem{AT2005}
{M.M.~Alexandrino and D.~T\"oben},
\newblock {\it Singular {R}iemannian foliations on simply connected spaces},
\newblock {\bf Differential Geom. Appl.}, to appear.

\bibitem{ALQ2001}
{R.~Ayala, F.F.~Lasheras and A.~Quintero},
\newblock {\it The equivariant category of proper $G$-spaces},
\newblock {\bf Rocky Mountain J. Math.}, 31:1111--1132, 2001.

\bibitem{Bartsch1993}
{T.~Bartsch},
\newblock   {\bf Topological methods for variational problems with symmetries},
\newblock {Lect. Notes in Math. Vol. 1560},  
\newblock {Springer--Verlag, Berlin}, 1993.

\bibitem{BH1984a}
{R.~Blumenthal and J.~Hebda},
\newblock {\it Ehresmann connections for foliations},
\newblock {\bf Indiana Univ. Math. J.}, 33:597--611, 1984.

\bibitem{BH1984b}
{R.~Blumenthal and J.~Hebda},
\newblock {\it Complementary distributions which preserve the leaf geometry and applications to totally geodesic foliations},
\newblock {\bf Quart. J. Math. Oxford Ser. (2)}, 35:383--392, 1984.

\bibitem{Bredon1972}
{G.~Bredon},
\newblock {\bf Introduction to compact transformation groups},
\newblock {Pure and Applied Mathematics, Vol. 46}, 
\newblock {Academic Press, New York}, 1972.

\bibitem{CN1985}
{C.~Camacho and A.~Lins~Neto},
\newblock {\bf Geometric Theory of Foliations},
\newblock {Translated from the Portuguese by Sue E. Goodman},
\newblock {Progress in Mathematics}, {Birkh\"auser Boston, MA}, 1985.

\bibitem{CandelConlon2000}
{A.~Candel and L.~Conlon},
\newblock {\bf Foliations I},
\newblock Amer. Math. Soc., Providence, RI, 2000.

\bibitem{Carriere1984}
{Y.~Carri{\`e}re},
\newblock {\it Flots riemanniens},
\newblock In {\bf Transversal structure of foliations (Toulouse, 1982)},
\newblock {Asterisque, 116, Soci\'et\'e Math\'ematique de France}, 1984, 31--52.

\bibitem{Colman1998}
{H.~Colman},
\newblock {\it Categor\'{\i}a LS en foliaciones},
\newblock {\bf Publicaciones del Departamento de Topolog\'{\i}a y Geometr\'{\i}a}, no. 93,
\newblock {Universidade de Santiago de Compostele}, 1998.

\bibitem{Colman2002a}
{H.~Colman},
\newblock {\it L{S}-categories for foliated manifolds},
\newblock In {\bf Foliations: Geometry and Dynamics (Warsaw, 2000)},
\newblock {World Scientific Publishing Co. Inc., River Edge, N.J.}, 2002:17--28.

\bibitem{Colman2002b}
{H.~Colman},
\newblock {\it Equivariant LS-category for finite group actions},
\newblock In {\bf Lusternik-Schnirelmann category and related topics (South Hadley, MA, 2001)},
\newblock {Contemp. Math. Vol. 316}, {Amer. Math. Soc., Providence, R.I.}, 2002, 35--40.

\bibitem{Colman2004}
{H.~Colman},
\newblock {\it Transverse {Lusternik--Schnirelmann} category of {Riemannian} foliations},
\newblock {\bf Topology Appl.}, 141:187--196, 2004.

\bibitem{Colman2006}
{H.~Colman},
\newblock {\it {Lusternik--Schnirelmann} category of {Orbifolds}},
\newblock {\bf preprint}, 2006.

\bibitem{CH2004}
{H.~Colman and S.~Hurder},
\newblock {\it L{S}-category of compact {H}ausdorff foliations},
\newblock {\bf Trans. Amer. Math. Soc.}, 356:1463--1487, 2004.

\bibitem{CM2001}
{H.~Colman and E.~Macias},
\newblock {\it Transverse Lusternik--Schnirelmann category of foliated manifolds},
\newblock {\bf Topology } Vol. {\bf 40} (2) (2001), 419-430.

\bibitem{CLOT2003}
{O.~Cornea, G.~Lupton, J.~Oprea, and D.~Tanr\'e},
\newblock {\bf Lusternik-Schnirelmann category},
\newblock {Mathematical Surveys and Monographs} {\bf 103}, American Mathematical Society, 2003.

\bibitem{Davis1978}
{M.W.~Davis},
\newblock {\it Smooth $G$-manifolds as collections of fiber bundles},
\newblock {\bf Pacific J. Math.}, 77:315--363, 1978.

\bibitem{DK2000}
{J.J.~Duistermaat and J.A.C.~Kolk},
\newblock {\bf Lie Groups},
\newblock {Universitext, Springer-Verlag, Berlin}, 2000.

\bibitem{Epstein1976}
{{D.B.A.}~Epstein},
\newblock {\it Foliations with all leaves compact},
\newblock {\bf Ann. Inst. Fourier (Grenoble)}, 26:265--282, 1976.

\bibitem{Epstein1984}
{{D. B. A.} Epstein},
\newblock {\it Transversely hyperbolic $1$-dimensional foliations},
\newblock In {\bf Transversal structure of foliations (Toulouse, 1982)},
\newblock {Asterisque, 116, Soci\'et\'e Math\'ematique de France}, 1984, 53--69.

\bibitem{Fadell1985}
{E.~Fadell},
\newblock {\it The equivariant Lusternik-Schnirelmann method for invariant functionals and relative cohomological index theories},
\newblock In {\bf Topological methods in nonlinear analysis}, {S\'em. Math. Sup. Vol. 95},  ed. A. Granas, 
\newblock {Presses Univ. Montr\'eal, Montreal, QC}, 1985.

\bibitem{FH1987}
{E.~Fadell and S.~Husseini},
\newblock {\it Relative cohomological index theories},
\newblock {\bf Adv. Math.}, 64:1--31, 1987.

\bibitem{FHT2001}
{Y.~Felix, S.~Halperin, and J.-C.~Thomas},
\newblock {\bf Rational homotopy theory},
\newblock {Graduate Texts in Mathematics}, Vol. 205,
\newblock {Springer-Verlag, new York}, 2001.

\bibitem{Ghys1984}
{{\'E.} Ghys},
\newblock {\it Feuilletages riemanniens sur les vari\'et\'es simplement connexes},
\newblock {\bf Ann. Inst. Fourier (Grenoble).}, 34:203--223, 1984.

\bibitem{Haefliger1971}
{A.~Haefliger},
\newblock {\it Homotopy and integrability},
\newblock In {\bf Manifolds--Amsterdam 1970 (Proc. Nuffic Summer School)},
\newblock {Lecture Notes in Mathematics, Vol. 197}, 
\newblock {Springer, Berlin}, 133--163, 1971.

\bibitem{Haefliger1985}
{A.~Haefliger},
\newblock {\it Pseudogroups of local isometries},
\newblock In {\bf Differential geometry (Santiago de Compostela, 1984)},
\newblock {Res. Notes in Math.}, Vol. 131:174--197,
\newblock {Pitman, Boston, MS}, 1985.

\bibitem{Haefliger1988}
{A.~Haefliger},
\newblock {\it Leaf closures in Riemannian foliations},
\newblock In {\bf A F\^ete of Topology},
\newblock {Academic Press, Boston, MA}, 1988, 3--32.

\bibitem{Haefliger1989}
{A.~Haefliger},
\newblock {\it Feuilletages riemanniens},
\newblock In {\bf S\'eminaire Bourbaki, Vol.\ 1988/89},
\newblock {Asterisque, 177-178, Soci\'et\'e Math\'ematique de France}, 1989, 183--197.

\bibitem{HS1990}
{A.~Haefliger and {\'E}.~Salem},
\newblock {\it Riemannian foliations on simply connected manifolds and actions of tori on orbifolds},
\newblock {\bf Illinois J. Math.}, 34:706--730, 1990.

\bibitem{HS1991}
{A.~Haefliger and {\'E}.~Salem},
\newblock {\it Actions of tori on orbifolds},
\newblock {\bf Ann. Global Anal. Geom.}, 9:37--59, 1991.

\bibitem{HsHs1967}
{W.-C.~Hsiang and W.-Y.~Hsiang},
\newblock {\it Differentiable actions of compact connected classical groups. I},
\newblock {\bf Amer. J. Math.}, 89:705--786, 1967.

\bibitem{Hurder1981}
{S.~Hurder},
\newblock {\it On the homotopy and cohomology of the classifying space of Riemannian   foliations},
\newblock {\bf Proc. Amer. Math. Soc.}, 81:485--489, 1981.

\bibitem{Hurder2006a}
{S.~Hurder},
\newblock {\it Category and compact leaves},
\newblock {\bf Topology Appl.}, 153:2135--2154, 2006.

\bibitem{HT2006b}
{S.~Hurder and D.~T\"{o}ben}, 
\newblock {\it The equivariant LS-category of polar actions},
\newblock {preprint}, 2006.

\bibitem{HT2006c}
{S.~Hurder and D.~T\"{o}ben},
\newblock {\it Residues and transverse LS category for Riemannian foliations},
\newblock {in preparation}, 2006.

\bibitem{HW2006}
{S.~Hurder and P.~Walczak},
\newblock {\it Compact foliations with finite transverse LS category},
\newblock {\bf Jour. Math. Soc. Japan}, to appear.

\bibitem{James1978}
{I.M.~James},
\newblock {\it On category, in the sense of Lusternik-Schnirelmann},
\newblock {\bf Topology} 17:331--348, 1978.

\bibitem{James1995}
{I.M.~James},
\newblock {\it Lusternik-Schnirelmann Category},
\newblock {\bf Chapter 27, Handbook of Algebraic Topology},1995,1293--1310.

\bibitem{Janich1968}
{K.~J{\"a}nich},
\newblock {\it On the classification of $O(n)$-manifolds},
\newblock {\bf Math. Ann.} 176:53--76, 1968.

\bibitem{KM1989}
{W.~ Krawcewicz  and W.~Marzantowicz},
\newblock {\it Lusternik-Schnirelman method for functionals invariant with respect to a finite group action},
\newblock {\bf J. Differential Equations}, 85:105--124, 1989.

\bibitem{LW2002}
{R.~Langevin and P.~Walczak}, 
\newblock {\it Transverse Lusternik-Schnirelmann category and non-proper leaves},
\newblock In {\bf Foliations: Geometry and Dynamics (Warsaw, 2000)}, 
\newblock {World Scientific Publishing Co. Inc., River Edge, N.J.}, 2002:351--354.

\bibitem{LP1976a}
{C.~Lazarov and J.~Paternack},
\newblock {\it Secondary characteristic classes for Riemannian foliations},
\newblock {\bf J. Differential Geometry} 11:365--385, 1976.

\bibitem{LP1976b}
{C.~Lazarov and J.~Paternack},
\newblock {\it Residues and characteristic classes for Riemannian foliations},
\newblock {\bf J. Differential Geometry} 11:599--612, 1976.

\bibitem{LS1934}
{L.~Lusternik and L.~Schnirelmann}, 
\newblock {\bf  M\'ethodes topologiques dans les Probl\`emes Variationnels}.
\newblock {Hermann, Paris}, 1934.

\bibitem{Marzantowicz1989}
{W.~Marzantowicz},
\newblock {\it A $G$-Lusternik-Schnirelman category of space with an action of a compact Lie group},
\newblock {\bf Topology}, 28:403--412, 1989.

\bibitem{Mei1983}
{X.-M.~Mei},
\newblock {\it Note on the residues of the singularities of a Riemannian foliation},
\newblock {\bf Proc. Amer. Math. Soc.} 89:359--366, 1983.

\bibitem{Millett1974}
{K.~Millett},
\newblock {\it Compact foliations},
\newblock In {\bf Differential topology and geometry (Proc. Colloq., Dijon, 1974)},
\newblock {Lect. Notes in Math. Vol. 484}, 277--287,1975. Springer--Verlag, New York and Berlin.

\bibitem{Molino1977}
{P.~Molino},
\newblock {\it \'Etude des feuilletages transversalement complets et applications},
\newblock {\bf Ann. Sci. \'Ecole Norm. Sup. (4)} 10:289--307, 1977.

\bibitem{Molino1982}
{P.~Molino},
\newblock {\it G\'eom\'etrie globale des feuilletages riemanniens},
\newblock {\bf Nederl. Akad. Wetensch. Indag. Math.} 44:45--76, 1982.

\bibitem{Molino1988}
{P.~Molino},
\newblock {\bf Riemannian foliations},
\newblock {Translated from the French by Grant Cairns, with appendices by Cairns, Y. Carri\`ere, \'E. Ghys, E. Salem and V. Sergiescu},
\newblock {Birkh\"auser Boston Inc., Boston, MA}, 1988.

\bibitem{Molino1994}
{P.~Molino},
\newblock {\it Orbit-like foliations},
\newblock In {\bf Geometric Study of Foliations, Tokyo 1993} (eds. Mizutani et al),
\newblock {World Scientific Publishing Co. Inc., River Edge, N.J.}, 1994, 97--119.

\bibitem{MM2003}
{I.~Moerdijk and J.~Mr{\v{c}}un},
\newblock {\bf Introduction to foliations and Lie groupoids},
\newblock {Cambridge Studies in Advanced Mathematics}, Vol. 91,
\newblock {Cambridge University Press, Cambridge, 2003}.

\bibitem{Palais1966b}
{R.S.~Palais}, 
\newblock {\it Lusternik-{S}chnirelman theory on {B}anach manifolds}, 
\newblock {\bf Topology},  5:115--132,  1966. 


\bibitem{PT1988}
{R.~Palais and C.-L.~Terng},
\newblock   {\bf Critical point theory and submanifold geometry},
\newblock {Lect. Notes in Math. Vol. 1353},  
\newblock {Springer--Verlag, Berlin}, 1988.

\bibitem{Richardson1998}
{K.~Richardson},
\newblock {\it The asymptotics of heat kernels on Riemannian foliations},
\newblock {\bf Geom. Funct. Anal.}, 8:356--401, 1998.

\bibitem{Richardson2001}
{K.~Richardson},
\newblock {\it The transverse geometry of $G$-manifolds and Riemannian foliations},
\newblock {\bf Illinois J. Math.}, 45:517--535, 2001.

\bibitem{Salem1988}
{{\'E}.~Salem},
\newblock {\it Une g\'en\'eralisation du th\'eor\`eme de Myers-Steenrod aux pseudogroupes d'isom\'etries},
\newblock {\bf Ann. Inst. Fourier (Grenoble)}, 38:185--200, 1988.

\bibitem{Satake1956}
{I.~Satake},
\newblock {\it On a generalization of the notion of manifold},
\newblock {\bf Proc. Nat. Acad. Sci. U.S.A. }, 42:359--363, 1956.

\bibitem{Stefan1974}
{P.~Stefan},
\newblock {\it Accessible sets, orbits, and foliations with singularities},
\newblock {\bf Proc. London Math. Soc. (3)}, 29:699-713, 1974.

\bibitem{Stefan1980}
{P.~Stefan},
\newblock {\it Integrability of systems of vector fields},
\newblock {\bf J. London Math. Soc. (2)}, 21:544--556, 1980.

\bibitem{tomDieck1987}
{T.~tom Dieck},
\newblock {\bf  Transformation groups},
\newblock {de Gruyter Studies in Mathematics, Vol. 8}, 1987.
 
\bibitem{Wolak2000}
{R.A.~Wolak},
\newblock {\it Basic forms for transversely integrable singular Riemannian foliations},
\newblock {\bf Proc. Amer. Math. Soc.}, 128:1543--1545, 2000.

\bibitem{Wolak2002}
{R.A.~Wolak},
\newblock {\it Critical leaves of basic functions for a singular Riemannian foliation},
\newblock {preprint 2002}.


\bibitem{Yamato1979}
{K.~Yamato},
\newblock {\it Sur la classe caract\'eristique exotique de Lazarov-Pasternack en codimension $2$},
\newblock {\bf C. R. Acad. Sci. Paris S\'er. A-B}, 289:A537--A540, 1979.

\bibitem{Yamato1981}
{K.~Yamato},
\newblock {\it Sur la classe caract\'eristique exotique de Lazarov-Pasternack en codimension $2$. {II}},
\newblock {\bf Japan. J. Math. (N.S.)}, 7:227--256, 1981.

\end{thebibliography}
\end{document}